\documentclass{article}

\usepackage{amsmath}   
\usepackage{amssymb}

\newtheorem{Th}{Theorem}[section] 
\newtheorem{Prop}{Proposition}[section]   

\newtheorem{Lemma}{Lemma}[section]   
\newtheorem{Coro}{Corollary}[section]

\newcommand{\finishproof}{\hfill $\Box$ \vspace{3mm}}
\newcommand{\R}{\mathbb{R}}
\newcommand{\Z}{\mathbb{Z}}
\newcommand{\C}{\mathbb{C}}
\newcommand{\T}{\mathbb{T}}

\newcommand{\F}{{\mathcal F}}
\newcommand{\W}{{\mathcal W}}

\newcommand{\V}{{\mathcal V}}
\newcommand{\D}{{\mathcal D}}
\newcommand{\id}{{\tt Id}}
\newcommand{\res}{\mathop{\tt Res}}

\newcommand{\spec}{\text{Spec}_{\text{per}}}
\newcommand{\LL}{L^2_c}
\newcommand{\LLb}{L^2_\bullet}
\newcommand{\LLr}{L^2_r}
\newcommand{\CC}{{\mathcal C}}
\newcommand{\al}{\alpha}
\newcommand{\ep}{\varepsilon}
\newcommand{\ga}{\gamma}
\newcommand{\la}{\lambda}
\newcommand{\si}{\sigma}
\newcommand{\tsi}{\tilde\sigma}
\newcommand{\de}{\delta}

\newcommand{\Ga}{\Gamma}
\newcommand{\De}{\Delta}

\begin{document}

\date{}

\title{On normalized differentials on families of curves of infinite genus}   
 
\author{T. Kappeler\footnote{Supported in part by the Swiss National Science Foundation.},
P. Lohrmann$^*$, and P. Topalov\footnote{Supported in part by NSF DMS-0901443}} 
  
\maketitle

\begin{abstract}  
\noindent We construct normalized differentials on families of curves of infinite genus.
Such curves are used to investigate integrable PDE's such as the focusing
Nonlinear Schr{\"o}dinger equation.
\end{abstract}

\section{Introduction}\label{sec:introduction}
In this paper we study a family of curves associated with the Zakharov-Shabat
operator ({\em ZS operator}),
\begin{equation}\label{e:ZS_intro}
L(\varphi):=i\,\left(
\begin{array}{cc}
1&0\\
0&-1
\end{array}
\right)\!\partial_x+
\left(
\begin{array}{cc}
0&\varphi_1\\
\varphi_2&0
\end{array}
\right)
\end{equation}
with periodic (or anti-periodic) boundary conditions.
We assume that $\varphi=(\varphi_1,\varphi_2)\in\LL$ where $\LL=L^2\times L^2$ and $L^2=L^2(\T,\C)$ is
the Hilbert space of square integrable complex-valued functions on the circle $\T:=\R/\Z$.
More specifically, we consider the curve,
\begin{equation}\label{e:Sigma_intro}
\CC_\varphi:=\{(\la,w)\in\C^2\,|\,w^2=\De(\la,\varphi)^2-4\}\,,
\end{equation}
where $\De(\la,\varphi)$ is the discriminant of the ZS operator. It is known (see e.g. \cite{LM} or \cite{GKP}) that
\[
\De(\la,\varphi)^2-4=-4\prod_{k\in\Z}\frac{(\la-\la_k^-(\varphi))(\la-\la_k^+(\varphi))}{\pi_k^2}
\]
where
\begin{equation}\label{e:pi_k}
\pi_k:=\left\{
\begin{array}{l}
k\pi,\,\,\,\,\,\,\mbox{if}\,\,\,k\ne 0\\
1,\,\,\,\,\,\,\,\,\,\,\mbox{if}\,\,\,k=0
\end{array}
\right.
\end{equation}
and where $(\la_k^\pm)_{k\in\Z}$ denotes the periodic eigenvalues of $L(\varphi)$, appropriately ordered and
listed with multiplicities -- see Section\,\ref{sec:ZS_operators} for details.
Note that the periodic spectrum $\spec L(\varphi)$ is pure point and hence consists only of eigenvalues.
We do {\em not} necessarily restrict our attention to potentials of real type,
\[
\varphi\in\LLr=\{(\varphi_1,\varphi_2)\in\LL\,|\,\varphi_1=\overline{\varphi_2}\}\,,
\]
for which $L(\varphi)$ is self-adjoint or
\[
\varphi\in i\LLr=\{(\varphi_1,\varphi_2)\in\LL\,|\,\varphi_1=-\overline{\varphi_2}\}\,.
\]
As a subset of $\C^2$, the curve $\CC_\varphi$ is a topological space whereas
\[
\CC_\varphi^\bullet:=\CC_\varphi\setminus\{(\la,0)\in\C^2\,|\,
\la\,\,\mbox{is a multiple periodic eigenvalue of}\,L(\varphi)\}
\]
is an open Riemann surface.

When $\varphi\in\LL$ varies the singularities of the curve $\CC_\varphi$ might change as well, i.e.,
additional multiple eigenvalues might emerge or a multiple eigenvalue could split up into several simple 
and/or multiple eigenvalues of smaller multiplicity. Our aim is to construct a family of holomorphic differentials
on $\CC_\varphi^\bullet$, depending analytically on $\varphi$ and normalized with respect to a properly chosen infinite
set of cycles on $\CC_\varphi^\bullet$. More specifically, we prove the following result.

Let $\varphi_*\in\LL$ and assume that the periodic spectrum of $L(\varphi_*)$ has only eigenvalues of algebraic
multiplicity one or two.
Then there exist an open neighborhood $\W$ of $\varphi_*$ in $\LL$ and a family of simple, closed, smooth, oriented
curves $\Ga_k\subseteq\C\setminus\spec L(\varphi_*)$, $k\in\Z$, so that the closures of the domains in $\C$,
bounded by the $\Ga_k$'s, are pairwise disjoint and for any $\varphi\in\W$ and $k\in\Z$
\begin{itemize}
\item[(H1)] the domain bounded by $\Ga_k$ contains precisely the two eigenvalues
$\la_k^\pm(\varphi)$;
\item[(H2)] there is a cycle $A_k$ on $\CC_\varphi^\bullet$ such that $\pi(A_k)=\Ga_k$ where $\pi$ is the projection
$\pi : \CC_\varphi^\bullet\to\C$, $(\la,w)\mapsto\la$.
\end{itemize}
A potential $\varphi\in\LL$ is called a {\em finite gap} potential if the periodic spectrum of $L(\varphi)$
has only finitely many simple eigenvalues. The set of finite gap potentials is dense in $\LL$ -- see e.g. \cite{KST2}.
We prove the following theorem:
\begin{Th}\hspace{-2mm}{\bf .}\label{Th:main}
Let $\varphi_*\in\LL$ be a finite gap potential such that all periodic eigenvalues of $L(\varphi_*)$ have algebraic
multiplicity at most two. Then there exist an open neighborhood $\W$ of $\varphi_*$ in $\LL$ and a family of
analytic functions $\zeta_n : \C\times\W\to\C$ so that for any $\varphi\in\W$ and $m,n\in\Z$,
\begin{equation}\label{e:normalization*}
\frac{1}{2\pi}\oint_{A_m}\frac{\zeta_n(\la,\varphi)}{\sqrt{\De(\la,\varphi)^2-4}}\,d\la=\de_{mn}\,.
\end{equation}
For any $n\in\Z$, the entire function $\zeta_n(\cdot,\varphi)$ has a product representation
\begin{equation}\label{e:zeta_product*}
\zeta_n(\la,\varphi)=-\frac{2}{\pi_n}\prod_{k\ne n}\frac{\tsi_k^{(n)}-\la}{\pi_k}\,,\,\,\,
\tsi_k^{(n)}=\tsi_k^{(n)}(\varphi),
\end{equation}
with $\tsi_k^{(n)}=k\pi+l^2(k)$ uniformly in $n\in\Z$ and $\varphi\in\W$, i.e.,
the sequences $(\tsi_k^{(n)}(\varphi)-k\pi)_{k \in \mathbb Z}$ are bounded in $l^2(\mathbb Z,\mathbb C)$
uniformly in $n\in\Z$ and $\varphi\in\W$.
Furthermore, there exists $N\ge 1$ such that for any $|k|\ge N+1$, $k\ne n$,
\begin{equation}\label{e:asymptotics*}
\tsi_k^{(n)}=\tau_k+O(|\la_k^+-\la_k^-|^2),\,\,\,\,\,\,\tau_k=\frac{\la_k^-+\la_k^+}{2},
\end{equation}
uniformly in $n\in\Z$ and uniformly in $\varphi\in\W$. For $|k|\le N$, $k\ne n$,
\[
\tsi_k^{(n)}\in\{\la\in\C\,|\,|\la|\le (N+1/4)\pi\}
\]
uniformly in $n\in\Z$ and uniformly in $\varphi\in\W$. Moreover, if $\la_k^-=\la_k^+$ and $k\ne n$ then
$\tau_k$ is a zero of the entire function $\zeta_n$.
\end{Th}

\vspace{0.2cm} 
 
\noindent\textit{Applications:} In \cite{KLT}, Theorem \ref{Th:main} is used to construct locally near generic potentials
action-angle coordinates for the focusing Nonlinear Schr{\"o}dinger equation significantly extending previous results
in this direction obtained in \cite{KLTZ} for the zero potential. See also \cite{AbMa} for the discussion of $1$-gap
and $2$-gap potentials. Such coordinates allow to obtain various results concerning well-posedness for
these equations and study their (Hamiltonian) perturbations-- see \cite{KT1,KT2} respectively \cite{KP}
where corresponding results for the KdV equation have been obtained.

\vspace{0.2cm} 

In the remainder of this paper we prove Theorem \ref{Th:main}.
In Section \ref{sec:ZS_operators}, we collect facts about the Zakharov-Shabat operator needed in the sequel.
In Section \ref{sec:finite_gaps}, we construct entire functions $\zeta_n(\cdot,\varphi)$ satisfying
\eqref{e:normalization*}-\eqref{e:asymptotics*} in the case $\varphi$ is a finite gap potential.
This construction leads to the analytic set-up for the proof of Theorem \ref{Th:main} explained
in Section \ref{sec:set-up}. Theorem \ref{Th:main} is then proved in Section \ref{sec:implicit_function} and
\ref{sec:uniformity} by the implicit function theorem.
By a similar approach a version of Theorem \ref{Th:main} has been proved in \cite{GKP} for potentials
in a (small) neighborhood of $\LLr$ in $\LL$. The approach in \cite{GKP} had to be significantly modified
as the non-degeneracy condition needed for applying the implicit function theorem is not satisfied for
the potentials considered in the present paper. 
To keep the paper relatively simple we decided to treat the extension of Theorem \ref{Th:main} to potentials which are
not necessarily finite gap potentials in a subsequent paper.
 
\vspace{0.2cm} 
 
\noindent\textit{Related results:} In the case of a Hill operator with all periodic eigenvalues simple,
hence with the corresponding curve $\CC_\varphi$ being a Riemann surfaces of infinite genus, the existence of normalized
holomorphic differentials $\omega_n$, $n\in\Z$, on $\CC_\varphi$ can be deduced from Hodge theory (see \cite{FKT,MSS}).
This construction does not work in the case of the ZS-operator as the differentials constructed in Theorem \ref{Th:main}
are not square integrable on $\CC_\varphi^\bullet$. Another obstacle for applying a similar approach to our situation
is related to the fact that the set of finite gap potentials of ZS-operators is dense in $\LL$ (\cite{KST2}), and hence any 
open neighborhood $\W$ of $\varphi_*$ in $\LL$ contains potentials with double eigenvalues. 
In addition, the analytic dependence of the normalized differentials on $\varphi\in\W$,
their precise form \eqref{e:zeta_product*} and the uniform localization of the zeroes \eqref{e:asymptotics*}
cannot be obtained from the general theory in \cite{AS}. Note also that a different approach was used in \cite{MV}
for potentials in $\LLr$.


\section{Zakharov-Shabat operator}\label{sec:ZS_operators}
Denote by $\LL$ the Cartesian product $L^2\times L^2$ where $L^2:=L^2(\T,\C)$ is the Hilbert space
of square integrable complex-valued functions on the circle $\T:=\R/\Z$. 
For $\varphi=(\varphi_1,\varphi_2)\in\LL$ consider
the Zakharov-Shabat operator,
\begin{equation}\label{e:ZS}
L(\varphi):=i\,\left(
\begin{array}{cc}
1&0\\
0&-1
\end{array}
\right)\!\partial_x+
\left(
\begin{array}{cc}
0&\varphi_1\\
\varphi_2&0
\end{array}
\right).
\end{equation}
For any $\lambda \in \mathbb C$, let $M=M(x,\lambda,\varphi)$ denote the fundamental $2\times 2$ matrix of $L(\varphi)$,
\begin{equation*}
L(\varphi)M=\lambda M,
\end{equation*} 
satisfying the initial condition $M(0,\lambda,\varphi)=\text{Id}_{2\times 2}$.

\vspace{0.2cm} 
 
\noindent\textit{Periodic spectrum:} Denote by $\spec(\varphi)$ the spectrum of $L(\varphi)$ with domain
\begin{equation*}
\text{dom}_{\text{per}}(L):=\lbrace F\in H^1_{loc}\times H^1_{loc}| \:  F(1)=\pm F(0)\rbrace
\end{equation*}
where $H^1_{loc}\equiv H^1_{loc}(\R,\C)$.
This spectrum coincides with the spectrum of $L(\varphi)$ considered on $[0,2]$ with
periodic boundary conditions. 

We say that two complex numbers $a,b$ are \textit{lexicographically ordered}, $a \preccurlyeq b$,
if $[\text{Re}(a) < \text{Re}(b)]$ or $[ \text{Re}(a)  = \text{Re}(b) \text{ and } \text{Im}(a) \leq \text{Im}(b)]$.
Similarly, $a\prec b$ if $a\preccurlyeq b$ and $a\ne b$.

The following propositions are well known\,--\,see e.g. \cite[Section 3]{GKP}, .
\begin{Prop}\hspace{-2mm}{\bf .}\label{Prop:counting_lemma}
For any $\varphi_*\in\LL$ there exist an open neighborhood $\W$ of $\varphi_*$ in $\LL$ and
an integer $N_0\ge 1$ such that for any $\varphi\in\W$, the following statements hold.
\begin{itemize}
\item[$(i)$] For any $k\in\Z$ with $|k|\ge N_0+1$, the disk $\{\la\in\C\,|\,|\la-k\pi|<\pi/4\}$ contains precisely two (counted with multiplicities)
periodic eigenvalues $\la_k^{-}\preccurlyeq\la_k^{+}$ of $L(\varphi)$.
\item[$(ii)$] The disk $\{\la\in\C\,|\,|\la|<(N_0+1/4)\pi\}$ contains precisely $4N_0+2$ periodic eigenvalues of $L(\varphi)$
counted with multiplicities.
\item[$(iii)$] There are no other periodic eigenvalues of $L(\varphi)$ than the ones listed in items $(i)$ and $(ii)$.
\end{itemize}
\end{Prop}
Let $\W$ be the neighborhood given by Proposition \ref{Prop:counting_lemma}.
Note that for any $\varphi\in\W$,
\[
...\preccurlyeq\la_{-k-1}^+\prec\la_{-k}^-\preccurlyeq\la_{-k}^+\prec...\,\,\,\mbox{and}\,\,\,
...\prec\la_k^-\preccurlyeq\la_k^+\prec\la_{k+1}^-\preccurlyeq...
\]
for any $k\ge N_0+1$. 
\begin{Prop}\hspace{-2mm}{\bf .}\label{Prop:periodic_spectrum}
For any $\varphi\in\W$, the periodic eigenvalues $(\la_k^\pm)_{|k|\ge N_0+1}$ of $L(\varphi)$ ordered as above satisfy
\begin{equation*}
\lambda_k^\pm(\varphi)=k\pi +l^2(k)
\end{equation*}
locally uniformly in $\varphi$, i.e.
$(\lambda_k^\pm(\varphi)-k\pi)_{|k|\ge N_0+1}$ is locally bounded in $l^2$. 
\end{Prop}

\vspace{0.2cm}

\noindent\textit{Ordering of eigenvalues:} 
Denote by $\LLb$ the set of potentials $\varphi\in\LL$ so that all periodic eigenvalues of $L(\varphi)$ are of algebraic multiplicity at most two.
Take $\varphi_*\in\LLb$ and choose $\W\subseteq\LLb$ and $N_0\ge 1$ as in Proposition \ref{Prop:counting_lemma}.
We keep the lexicographic order of the eigenvalues that appear in item $(i)$. 
The remaining $4N_0+2$ eigenvalues in item $(ii)$ are grouped in lexicographically ordered pairs of two,
$\la_k^-\preccurlyeq\la_k^+$, $|k|\le N_0$, in the following way: Choose $\la_{-N_0}^-$ to be
the smallest\footnote{With respect to the partial order $\preccurlyeq$ introduced above.} eigenvalue bigger than $\la_{-N_0-1}^+$.
If $\la_{-N_0}^-$ is double, then set $\la_{-N_0}^+:=\la_{-N_0}^-$, otherwise
denote by $\la_{-N_0}^+$ the smallest {\em simple} eigenvalue which is bigger than $\la_{-N_0}^-$.
Next, define $\la_{-N_0+1}^-$ to be the smallest eigenvalue bigger than $\la_{-N_0}^-$ and different from 
$\la_{-N_0}^+$ and determine $\la_{-N_0+1}^+$ in the same fashion as $\la_{-N_0}^+$.
Continuing in this way we arrive at a listing of the $4N_0+2$ eigenvalues so that $(\la_k^-)_{|k|\le N_0}$ are
in strictly increasing order,
\[
\la_{-N_0}^-\prec\la_{-N_0+1}^-\prec...\prec\la_{N_0}^-
\]
and so that all double eigenvalues form a pair.
By shrinking the neighborhood $\W$ if necessary, we choose for any $k\in\Z$ a simple counterclockwise oriented
$C^1$-smooth closed curve $\Ga_k$ in $\C$ such that the closures $\D_k$ of the domains bounded by the $\Ga_k$'s are
pairwise disjoint and such that for any $\varphi\in\W$ and $k\in\Z$, the domain bounded by $\Ga_k$ contains
precisely the two eigenvalues $\la_k^-(\varphi)\preccurlyeq\la_k^+(\varphi)$ as well as
a continuously differentiable simple curve $G_k=G_k(\varphi)$ connecting $\la_k^-(\varphi)$ with $\la_k^+(\varphi)$.
In the case $\la_k^-(\varphi)=\la_k^+(\varphi)$, $G_k$ is chosen to be the constant curve $\la_k^-(\varphi)$.
For any $|k|\ge N_0+1$, we choose $\Ga_k$ to be the counterclockwise oriented boundary of the disk
\[
D_k:=\{\la\in\C\,|\,|\la-k\pi|\le\pi/4\}
\]
and $G_k:=\{(1-t)\la_k^-(\varphi)+t\la_k^+(\varphi)\,|\,t\in[0,1]\}$ whereas for $|k|\le N_0$, $\Ga_k$
is chosen to be contained in the disk $\{\la\in\C\,|\,|\la|\le(N_0+\frac{1}{4})\pi\}$. 
Define
\[
\tau_k:=(\la_k^-+\la_k^+)/2\,\,\,\,\mbox{and}\,\,\,\,\ga_k:=\la_k^+-\la_k^-\,.
\]

\vspace{0.2cm}

\noindent\textit{Discriminant:} 
Let $\Delta(\lambda, \varphi) := \mathop{\tt tr}M(1,\la,\varphi)$ be the trace of $M(1,\la,\varphi)$.
It is well known that $\Delta(\la,\varphi)$ is an analytic function on
$\mathbb{C}\times\LL$ (cf. \cite[Section 3]{GKP}).
The proof of the following statement can be find in \cite[Proposition 3.4]{GKP}. 
\begin{Prop}\hspace{-2mm}{\bf .}\label{Prop:discriminant}
For any $\varphi\in\W$ and any $\lambda \in \mathbb C$,
\begin{equation*} \Delta^2(\la,\varphi)-4=
-4\prod_{k\in\Z}\frac{\left(\la-\la_k^+(\varphi)\right)\left(\la-\la_k^-(\varphi)\right)}{\pi_k^2}.
\end{equation*} 
\end{Prop} 

\vspace{0.3cm}

\noindent\textit{Standard\,\&\,canonical roots:} Let $a,b\in\C$.
Denote by $\sqrt[+]{z}$ the principal branch of the square root defined on
$\C\setminus\{z\in\R\,|\,z\le 0\}$ by $\sqrt[+]{1}=1$.
We define the standard root of $(\la-a)(\la-b)$ by the following relation
\begin{equation}\label{e:standard_root1}
\sqrt[s]{(\la-a)(\la-b)}=-\la\sqrt[+]{\Big(1-\frac{a}{\la}\Big)\Big(1-\frac{b}{\la}\Big)}
\end{equation}
for all $\la\in\C\setminus\{0\}$ such that $\left|\frac{a}{\la}\right|\le 1/2$ and $\left|\frac{b}{\la}\right|\le 1/2$.
Let $G_{[a,b]}$ be an arbitrary continuous simple curve connecting $a$ and $b$.
By analytic extension, \eqref{e:standard_root1} uniquely defines a holomorphic function on $\C\setminus G_{[a,b]}$,
that we call the {\em standard root} of $(\la-a)(\la-b)$ on $\C\setminus G_{[a,b]}$.
One has the asymptotic formula
\begin{equation}\label{e:standard_root2}
\sqrt[s]{(\la-a)(\la-b)}\,\,\sim\,-\la,\,\,\,\,\mbox{as}\,\,\,\,\,\,\,|\la|\to\infty.
\end{equation}
For any $\varphi\in\W$ and $\la\in\C\setminus\sqcup_{k\in\Z}G_k$ we define the {\em canonical root} of $\Delta(\la,\varphi)^2-4$
\begin{equation}\label{e:canonical_root}
\sqrt[c]{\Delta(\la,\varphi)^2-4}:=2i\prod_{k\in\Z}
\frac{\sqrt[s]{\left(\la-\la_k^+(\varphi)\right)\left(\la-\la_k^-(\varphi)\right)}}{\pi_k}\,.
\end{equation}

The proof of the following lemma is straightforward and hence omitted.
\begin{Lemma}\hspace{-2mm}{\bf .}\label{Lem:canonical_root}
For any $\varphi\in\W$, the canonical root \eqref{e:canonical_root} defines a holomorphic function on
$\C\setminus(\sqcup_{k\in\Z}G_k)$.
\end{Lemma}
For any $\varphi\in\W$, consider the curve
\begin{equation}
\CC_\varphi:=\{(\la,w)\in\C^2\,|\,w^2=\Delta(\la,\varphi)^2-4\}
\end{equation}
as well as its canonical branch
\begin{equation}
\CC_\varphi^c:=\{(\la,w)\in\C^2\,|\,\la\in\C\setminus(\sqcup_{k\in\Z}G_k), w=\sqrt[c]{\De(\la,\varphi)^2-4}\}\,.
\end{equation}
The cycle $A_m$, $m\in\Z$, introduced in Section \ref{sec:introduction}, is now defined more precisely as
the cycle on $\CC_\varphi^c$ whose projection onto $\C$ is $\Ga_m$, $A_m=\pi^{-1}(\Ga_m)\cap\CC_\varphi^c$, and
$\pi : \CC_\varphi\to\C$, $(\la,w)\mapsto\la$, denotes the projection onto $\C$.

\section{Finite gap potentials}\label{sec:finite_gaps}
In this section we construct entire functions $\zeta_n(\cdot,\varphi)$ satisfying
\eqref{e:normalization*}-\eqref{e:zeta_product*} for an arbitrary finite gap potentials in $\LLb$.
The treatment of this special case will lead to the set-up of the proof of Theorem \ref{Th:main}, discussed
in Section \ref{sec:set-up}.
To state the result more precisely, introduce for $\varphi\in\LLb$
\[
J\equiv J(\varphi):=\{n\in\Z\,|\,\la_n^-(\varphi),\,\la_n^+(\varphi)\,\,\mbox{are simple}\}\,.
\]
Then, for any $n\in\Z\setminus J$, $\la_n^-(\varphi)=\la_n^+(\varphi)$ is a double eigenvalue of $L(\varphi)$.
By definition, $\varphi\in\LLb$ is a finite gap potential if $J(\varphi)$ is finite.
\begin{Th}\hspace{-2mm}{\bf .}\label{Th:main_finite_gap}
Let $\varphi\in\LLb$ be a finite gap potential. Then for any $n\in\Z$, there exists an entire function
$\zeta_n(\cdot,\varphi)$ so that for any $m\in\Z$,
\begin{equation}\label{e:normalization_relations_fg}
\frac{1}{2\pi}\oint_{A_m}\frac{\zeta_n(\la,\varphi)}{\sqrt{\De(\la,\varphi)^2-4}}\,d\la=\de_{mn}
\end{equation}
where $A_m$, $m\in\Z$, are the cycles on the canonical sheet $\CC_\varphi^c$ of $\CC_\varphi$ introduced at
the end of Section \ref{sec:ZS_operators}.
Further, $\zeta_n(\la,\varphi)$ is of the form
\begin{equation}
\zeta_n(\la,\varphi)=
-\frac{2}{\pi_n}\Big(\prod_{j\in J,j\ne n}\pi_j\Big)^{-1}P_n(\la,\varphi)\prod_{j\notin J,j\ne n}\frac{\tau_j-\la}{\pi_j}
\end{equation}
where $P_n(\la,\varphi)$ is a polynomial in $\la$ of degree $d$,
$P_n(\la,\varphi)=(-\la)^{d}+...$, with $d=|J|$ if $n\notin J$ and $d=|J|-1$ otherwise.
\end{Th}
We begin with a few preparations for the proof of Theorem \ref{Th:main_finite_gap}.

For $\varphi\in\LLb$, $(\la_k^\pm\equiv\la_k^\pm(\varphi))_{k\in\Z}$, and $J=J(\varphi)$ as above,
introduce the compact Riemann surface 
\begin{equation}\label{e:Sigma_J}
\Sigma_J:=\CC_J\sqcup\{\infty^\pm\}
\end{equation}
of genus $|J|-1$ that is obtained by compactifying the affine curve
\begin{equation}\label{e:C_J}
{\mathcal C}_J:=\{(\la,w)\in\C^2\,|\,w^2=\prod_{k\in J}(\la-\la_{k}^{-})(\la-\la_{k}^{+})\}
\end{equation}
in a standard way by adding two points $\infty^{\pm}$ at infinity -- one for each sheet.
The charts in open neighborhoods of $\infty^{\pm}$ are defined by the local parameter $z=1/\la$.
Let $A_m=A_m^J$, $m\in\Z$, be cycles on the canonical sheet of $\Sigma_J$,
\begin{equation}\label{e:Sigma_J^c}
\Sigma_J^c:=\{(\la,w)\in\C^2\,|\,\la\in\C\setminus(\sqcup_{k\in J}G_k),
w=\prod_{k\in J}\sqrt[s]{(\la-\la_k^{-})(\la-\la_k^{+})}\}\sqcup\{\infty^+\}\,,
\end{equation}
such that $\pi_J(A_m)=\Ga_m$, where $\pi_J$ is the projection onto the first component
$\pi_J : {\mathcal C}_J\to\C$, $(w,\la)\mapsto\la$. Consider the differential on $\CC_J$, 
\begin{equation}\label{e:ansatz1}
\chi:=-\frac{(-\la)^{l-1}}{i\,\sqrt{\prod_{k\in J}(\la-\la_k^{-})(\la-\la_k^{+})}}\,d\la
\end{equation}
where $l:=|J|$.
The proof of the following lemma is straightforward and hence is omitted.
\begin{Lemma}\hspace{-2mm}{\bf .}\label{Lem:infinity1}
The differential $\chi$ extends to a meromorphic differential on $\Sigma_J$ with precisely two poles.
They are situated at $\infty^\pm$ and have residues $\res\limits_{\infty^\pm}\chi=\pm i$.
\end{Lemma}
For any $n\in J$, consider the basis of holomorphic differentials on $\Sigma_J$,
\begin{equation}\label{e:om_km}
\omega_{ns}:=\frac{P_{ns}(\la)}{\sqrt{\prod_{k\in J}(\la-\la_k^{-})(\la-\la_k^{+})}}\,d\la\,,\,\,\,\,\,\,\,
s\in J\setminus\{n\},
\end{equation}
where $P_{ns}(\la)$ are polynomials of degree $\deg P_{ns}\le l-2$, normalized by the conditions,
\begin{equation}\label{e:finite_normalization}
\oint_{A_m}\omega_{ns}=\de_{ms}\,,\,\,\,\,\,\,\,\,\, m,s\in J\setminus\{n\}\,. 
\end{equation}
For any $n\in J$, introduce the differential,
\begin{equation}\label{e:chi_k}
\chi_n:=\chi-\!\!\!\!\!\!\sum_{s\in J\setminus\{n\}}\!\!\!\!c_s\omega_{ns}
\end{equation}
where $c_s:=\oint_{A_s}\chi$. Note that $\chi_n$ has the same poles as $\chi$ and their residues coincide with the ones of
$\chi$. Furthermore, any of the cycles $A_k$, $k\in J$, is homologous to a connected component of the boundary of
$\Sigma_J^c$ in $\Sigma_J$.
By Lemma \ref{Lem:infinity1}, the closure of $\Sigma_J^c$ in $\Sigma_J$ contains precisely one pole of $\chi_n$.
It is situated at $\infty^+$ and is of order one with residue $i$.
For $m\in\Z\setminus J$ note that the cycle $A_m\equiv A_m^J$ bounds a disk which is contained in
$\Sigma_J^c\setminus\infty^+$.
As $\chi_n$ is a holomorphic 1-form on $\Sigma_J^c\setminus\infty^+$ one concludes that $\oint_{A_m}\chi_n=0$
for those $m$'s.
Furthermore, for any $m\in J\setminus\{n\}$ it follows from \eqref{e:chi_k} that again $\oint_{A_m}\chi_n=0$.
For $m=n$, one then concludes from Stokes' formula,
\[
\sum_{m\in J}\oint_{A_m}\chi_n+2\pi\,i\res_{\infty^+}\chi_n=0\,,
\]
that $\oint_{A_n}\chi_n=2\pi$. Summarizing, we get that for any $n\in J$ and for any $m\in\Z$,
\begin{equation}\label{e:infinite_normalization}
\frac{1}{2\pi}\oint_{A_m}\chi_n=\de_{mn}\,.
\end{equation}
By construction,
\begin{equation}\label{e:chi_k-summed}
\chi_n=-\frac{(-\la)^{l-1}+\al^{(n)}_1\la^{l-2}+...+\al^{(n)}_{l-1}}
{i\,\sqrt{\prod_{k\in J}(\la-\la_k^{-})(\la-\la_k^{+})}}\,d\la
\end{equation}
where $\al^{(n)}_j$ ($1\le j\le l-1$) are complex numbers.

\vspace{0.3cm}

To continue, let us consider the case where $n\in\Z\setminus J$. Then $\la_n^-=\la_n^+=\tau_n$ and we introduce
the following differential on $\CC_J$,
\begin{equation}\label{e:ansatz2}
{\tilde\chi}_n:=
\frac{(-\la)^l+\varepsilon_n(\la+e_n)^{l-1}}
{(\la-\tau_n)\,i\,\sqrt{\prod_{k\in J}(\la-\la_k^{-})(\la-\la_k^{+})}}\,d\la\,,
\end{equation}
where 
\[
e_n=\left\{
\begin{array}{l}
0,\,\,\,\tau_n\ne 0,\\
1,\,\,\,\tau_n=0,
\end{array}
\right.
\]
and $\varepsilon_n$ is a complex number chosen so that 
\begin{equation}\label{e:chi_normalization}
\frac{(-\tau_n)^l+\varepsilon_n(\tau_n+e_n)^{l-1}}
{i\,\prod_{k\in J}\sqrt[s]{(\tau_n-\la_k^{-})(\tau_n-\la_k^{+})}}=-i\,.
\end{equation}
Note that with the above definition of $e_n$, $\varepsilon_n$ is well defined by \eqref{e:chi_normalization}.
Let $\tau_n^\pm$ be the two points on $\Sigma_J$ so that $\pi_J(\tau_n^\pm)=\tau_n$ with
$\tau_n^+$ lying on the canonical branch $\Sigma_J^c$.
The following lemma follows easily from the normalization \eqref{e:chi_normalization}.
\begin{Lemma}\hspace{-2mm}{\bf .}\label{Lem:infinity2} For any $n\in\Z\setminus J$, the differential ${\tilde\chi}_n$
extends to a meromorphic differential on $\Sigma_J$ with precisely four poles.
They are located at $\tau_n^{\pm}$ and $\infty^\pm$ and their residues are $\res\limits_{\infty^\pm}{\tilde\chi}_n=\pm i$ and
$\res\limits_{\tau_n^{\pm}}{\tilde\chi}_n=\mp i$.
\end{Lemma}
For any $n\in\Z\setminus J$ define the differential,
\begin{equation}\label{e:chi_n}
\chi_n:={\tilde\chi}_n-\!\!\!\!\!\!\sum_{m\in J\setminus\{k_*\}}\!\!\!\!c_m^n\omega_{{k_*}m}\,,\,\,\,\,\,\,\,\,
c_m^n:=\oint_{A_m}{\tilde\chi}_n
\end{equation}
where $k_*$ is a fixed, but arbitrary number in $J$.
Recall that for any $k\in J$, the cycle $A_k$ is homologous to a connected component of the boundary of
$\Sigma_J^c$ in $\Sigma_J$. By Lemma \ref{Lem:infinity2}, the closure of $\Sigma_J^c$ in $\Sigma_J$ contains
precisely two poles of $\chi_n$, at $\tau_n^+$ and $\infty^+$, with residues $-i$ and $i$ respectively.
Arguing as above -- in particular using again Stokes' formula -- we conclude that $\oint_{A_k}\chi_n=0$ for any $k\in J$.
Using \eqref{e:chi_normalization} we obtain that for any $m\in\Z$,
\begin{equation}\label{e:infinite_normalization2}
\frac{1}{2\pi}\oint_{A_m}\chi_n=\de_{mn}\,.
\end{equation}
By construction,
\begin{equation}\label{e:chi_n-summed}
\chi_n=\frac{(-\la)^l+\al^{(n)}_1\la^{l-1}+...+\al^{(n)}_l}
{(\la-\tau_n)\,i\,\sqrt{\prod_{k\in J}(\la-\la_k^{-})(\la-\la_k^{+})}}\,d\la
\end{equation}
where $\al^{(n)}_j$ ($1\le j\le l$) are complex numbers.

\vspace{0.4cm}

\noindent{\em Proof of Theorem \ref{Th:main_finite_gap}.}
If $n\in J$ define 
\begin{equation}\label{e:defP_n-finite}
P_n(\la,\varphi):=(-\la)^{l-1}+\al_1^{(n)}\la^{l-2}+...+\al^{(n)}_{l-1}
\end{equation}
with $\al_j^{(n)}$ ($j=1,...,l-1$) as in \eqref{e:chi_k-summed}.
Then for any $m\in\Z$, the entire function
\[
\zeta_n(\la,\varphi)=
-\frac{2}{\pi_n}\Big(\prod_{j\in J,j\ne n}\pi_j\Big)^{-1}P_n(\la,\varphi)\prod_{j\notin J}\frac{\tau_j-\la}{\pi_j}
\]
satisfies the normalizing relation \eqref{e:normalization_relations_fg}.
Similarly, if $n\in\Z\setminus J$, define
\begin{equation}\label{e:defP_n-infinite}
P_n(\la,\varphi):=(-\la)^l+\al_1^{(n)}\la^{l-1}+...+\al^{(n)}_l
\end{equation}
with $\al_j^{(n)}$ ($j=1,...,l$) as in \eqref{e:chi_n-summed}.
Again, for any $m\in\Z$,
\[
\zeta_n(\la,\varphi)=
-\frac{2}{\pi_n}\Big(\prod_{j\in J}\pi_j\Big)^{-1}P_n(\la,\varphi)\prod_{j\notin J,j\ne n}\frac{\tau_j-\la}{\pi_j}
\]
satisfies the normalisation condition \eqref{e:normalization_relations_fg}.
\finishproof

\vspace{0.3cm}

\begin{Prop}\hspace{-2mm}{\bf .}\label{Prop:a-bounded}
There exists $T>0$ so that the coefficients $\al_j^{(n)}$, $1\le j\le l$, of the polynomial $P_n(\la,\varphi)$
in \eqref{e:defP_n-infinite} satisfy 
\[
\sup_{n\in\Z\setminus J,\,1\le j\le l}|\al_j^{(n)}|\le T\,.
\]
\end{Prop}
{\em Proof.} Let $n\in\Z\setminus J$.
It follows from the normalization condition \eqref{e:chi_normalization},
the asymptotic formula $\tau_n=n\pi+o(1)$ (see Proposition \ref{Prop:periodic_spectrum}),
and the property \eqref{e:standard_root2} of the $s$-root that $\varepsilon_n$,
defined in \eqref{e:chi_normalization}, satisfies for $n\to\pm\infty$
\begin{eqnarray}
\varepsilon_n&=&\Big(\prod_{k\in J}\sqrt[s]{(\tau_n-\la_k^{-})(\tau_n-\la_k^{+})}-(-\tau_n)^l\Big)/\tau_n^{l-1}
\nonumber\\
&=&(-1)^l\Big(\tau_n\Big(1+O(1/n)\Big)-\tau_n\Big)=O(1)\,.\label{e:eps-bounded}
\end{eqnarray}
This estimate together with \eqref{e:ansatz2} and the second formula in \eqref{e:chi_n} imply that
\begin{equation}\label{e:a_n-bounded}
c_m^n=O(1/n)
\end{equation}
uniformly in $m\in J\setminus\{k_*\}$. 
Furthermore, by the definition \eqref{e:chi_n} of ${\tilde\chi}_n$ and \eqref{e:om_km},
\[
{\tilde\chi}_n=\frac{(-\la)^l+\varepsilon_n(\la+e_n)^{l-1}+(\la-\tau_n)\,i\sum_{m\in J\setminus\{k_*\}}c_m^n P_{k_*m}(\la)}
{(\la-\tau_n)\,i\,\sqrt{\prod_{k\in J}(\la-\la_k^{-})(\la-\la_k^{+})}}\,.
\]
Recall that the polynomials $P_{k_*m}(\la)$ are of degree $\le l-2$ and their coefficients are independent of $n\in\Z$.
Hence, by \eqref{e:eps-bounded}, \eqref{e:a_n-bounded}, and $\tau_n=n\pi+o(1)$ as $n\to\infty$,
the polynomials defined by \eqref{e:defP_n-infinite},
\[
P_n(\la)=(-\la)^l+\varepsilon_n(\la+e_n)^{l-1}+\!\!\!\!\sum_{m\in J\setminus\{k_*\}}\!\!\!c_m^n \la i P_{k_*m}(\la)-
\!\!\!\!\sum_{m\in J\setminus\{k_*\}}\!\!\!c_m^n \tau_n i P_{k_*m}(\la)
\]
have coefficients bounded uniformly in $n\in\Z\setminus J$.
\finishproof

Let $\varphi_*\in\LLb$ be a finite gap potential and let $J\equiv J(\varphi_*)\subseteq\Z$ be the finite subset of indices so that
$\la_k^-\prec\la_k^+$, $k\in J$, are the simple periodic eigenvalues of $L(\varphi_*)$.
Choose an open neighborhood $\W$ of $\varphi_*$ in $\LLb$, $N_0\ge 1$ and cycles $\Ga_m$ and $A_m$, $m\in\Z$,
as in Section \ref{sec:ZS_operators}.
If necessary, choose $N_0\ge 1$ larger so that the disk $\{\la\in\C\,|\,|\la|\le(N_0+\frac{1}{4})\pi\}$ contains all the simple eigenvalues
$\{\la_k^\pm\,|\,k\in J\}$ of $L(\varphi_*)$. It follows from Theorem \ref{Th:main_finite_gap} that for any $N\ge N_0$
and $n\in\Z$, the entire function $\zeta_n(\cdot,\varphi_*)$ can be written as follows
\begin{equation}\label{e:ansatz_finite_gap}
\zeta_n(\la,\varphi_*)=-\frac{2}{\pi_n}\Big(\prod_{|j|\le N,j\ne n}\pi_j\Big)^{-1}P_n^N(\la,\varphi_*) \!\!\!\!  
\prod_{|j|\ge N+1,j\ne n}\!\!\!\frac{\tau_j-\la}{\pi_j}
\end{equation}
where
\begin{equation}\label{e:P_n^N}
P_n^N(\la,\varphi_*):=P_n(\la,\varphi_*)\!\!\!\!\!\!\prod_{|j|\le N,j\notin J\cup\{n\}}\!\!\!(\tau_j-\la),\,\,\,
\,\,\,\,\,\,\,\,\,\,\,\tau_j:=\tau_j(\varphi_*)\,.
\end{equation}
Proposition \ref{Prop:a-bounded} implies the following corollary.
\begin{Coro}\hspace{-2mm}{\bf .}\label{Coro:K}
Let $\varphi_*\in\LLb$, $J\equiv J(\varphi_*)\subseteq\Z$, and $N_0\ge 1$ be as above. Then there exist $N\ge N_0$ and a compact set 
${\mathcal K}\subseteq\C^{2N+1}$ so that the following statements hold:
\begin{itemize}
\item[$(i)$] for any $n\in\Z$ and $|j|\ge N+1$, the double eigenvalue $\tau_j$ is {\em not} a zero of the polynomial
$P_n^N(\la,\varphi_*)$;
\item[$(ii)$] for any $|n|\ge N+1$, the coefficient vector $(a^{(n)}_j)_{1\le j\le 2N+1}$ of the polynomial 
$P_n^N(\la,\varphi_*)=(-\la)^{2N+1}+a_1^{(n)}\la^{2N}+...+a_{2N+1}^{(n)}$ of \eqref{e:P_n^N},
lies in ${\mathcal K}$;
\item[$(iii)$] for any $a=(a_j)_{1\le j\le 2N+1}$ in ${\mathcal K}$, the zeroes of the polynomial
\[
Q^N(\la,a):=(-\la)^{2N+1}+a_1\la^{2N}+...+a_{2N+1}
\]
are contained in the disk $\{\la\in\C\,|\,|\la|\le(N+\frac{1}{4})\pi\}$.
In particular, for any $a\in{\mathcal K}$ and any $|j|\ge N+1$, the double eigenvalue $\tau_j$ is
{\em not} a zero of the polynomial $Q^N(\la,a)$.
\end{itemize}
\end{Coro}
\noindent{\em Proof.}
By Proposition \ref{Prop:a-bounded}, there exists $T>0$ so that, for any $n\in\Z\setminus J$, the coefficient vector
$\al^{(n)}=(\al^{(n)}_j)_{1\le j\le l}$ of the polynomial $P_n(\la)\equiv P_n(\la,\varphi_*)$ is bounded,
\[
\max_{1\le j\le l}|\al^{(n)}_j|\le T\,.
\]
Choose $N\ge N_0$ so that for any $(\al_j)_{1\le j\le l}\in\C^l$ with $\max\limits_{1\le j\le l}|\al_j|\le T$ the zeroes of
the polynomial $Q_\al(\la):=(-\la)^l+\al_1\la^{l-1}+...+\al_l$ are contained in the disk 
$B_N=\{\la\in\C\,|\,|\la|\le(N+\frac{1}{4})\pi\}$.
By choosing $N$ larger if necessary we may assume that for any $n\in J$,
the zeroes of the polynomial $P_n(\la)\equiv P_n(\la,\varphi_*)$ are contained in $B_N$ as well.
In view of the definition \eqref{e:P_n^N} of $P_n^N$, this choice of $N\ge N_0$ implies $(i)$.
To define the set ${\mathcal K}$, introduce for any vector $\al=(\al_j)_{1\le j\le l}$ the polynomial
\[
Q_\al(\la)\cdot\!\!\!\!\!\!\!\prod_{|j|\le N,j\notin J}\!\!\!\!(\tau_j-\la)
=(-\la)^{2N+1}+a_1\la^{2N}+...+a_{2N+1}\,.
\]
Then define the coefficient map
\[
\Phi : \C^l\to\C^{2N+1},\,\,\,\al\mapsto a
\]
where $a=(a_j)_{1\le j\le 2N+1}$ and let ${\mathcal K}$ be the image of 
${\mathcal B}_T^l=\{\al\in\C^l\,|\,\max\limits_{1\le j\le l}|\al_j|\le T\}$ by the map $\Phi$,
\[
{\mathcal K}:=\Phi({\mathcal B}_T^l)\subseteq\C^{2N+1}\,.
\]
As ${\mathcal B}_T^l\subseteq\C^l$ is compact and $\Phi$ is continuous it follows that ${\mathcal K}$ is compact.
By construction, item $(ii)$ and $(iii)$ hold.
\finishproof

\section{Analytic set-up}\label{sec:set-up}
Let $\varphi_*\in\LLb$ be a finite gap potential.
As in Section\,\,\ref{sec:ZS_operators} choose an open neighborhood $\W$ of $\varphi_*$ in $\LL$,
$N_0\ge 1$, and, for any $k\in\Z$, the cycle $\Ga_k$, as well as the curve $G_k\equiv G_k(\varphi)$
connecting the pair of eigenvalues $\la_k^-(\varphi)\preccurlyeq\la_k^+(\varphi)$, $\varphi\in\W$.
We want to construct analytic functions $\zeta_n(\la,\varphi)$, $n\in\Z$, on $\C\times\W$
with the neighborhood $\W$ shrinked if necessary such that for any $m,n\in\Z$,
\begin{equation}
\frac{1}{2\pi}\oint_{\Ga_m}\frac{\zeta_n(\la,\varphi)}{\sqrt[c]{\De(\la,\varphi)^2-4}}\,d\la=\de_{mn}\,.
\end{equation}
In order to do this we make an ansatz for $\zeta_n(\la,\varphi)$ and then determine the parameters involved
by applying the implicit function theorem. The ansatz is suggested by \eqref{e:ansatz_finite_gap} and \eqref{e:P_n^N}.

\vspace{0.2cm}

\noindent\textit{Ansatz:} Let $N\ge N_0$ be a given integer.
For any $|n|\le N$, define the entire function on $\C\times l^2\times\C^{2N}$,
\begin{equation}\label{e:f_n1}
f_n(\la,\si,a):=-\frac{2}{\pi_n}
\Big(\prod_{|j|\le N,j\ne n}\pi_j\Big)^{-1}Q^N_n(\la,a)\prod_{|j|\ge N+1}\frac{\tsi_j-\la}{\pi_j}\,,
\end{equation}
where $\si:=(\si_j)_{|j|\ge N+1}\in l^2$, $\tsi_j:=j\pi+\si_j$, and
\[
Q^N_n(\la,a):=(-\la)^{2N}+a_1\la^{2N-1}+...+a_{2N}\,,
\]
with $a:=(a_1,...,a_{2N})$ in $\C^{2N}$. 
Here and in the sequel, $l^2$ denotes the Hilbert space of complex valued sequences $(x_k)_{k\in I}$ with index set
$I\subseteq\Z$. It will be clear from the context what $I$ is.
Similarly, for $|n|\ge N+1$, define the entire function on $\C\times l^2\times\C^{2N+1}$,
\begin{equation}\label{e:f_n2}
f_n(\la,\si,a):=-\frac{2}{\pi_n}\Big(\prod_{|j|\le N}\pi_j\Big)^{-1}Q^N_n(\la,a)\prod_{|j|\ge N+1, j\ne n}\frac{\tsi_j-\la}{\pi_j}\,,
\end{equation}
where  $\si:=(\si_j)_{|j|\ge N+1,j\ne n}\in l^2$, $\tsi_j:=j\pi+\si_j$, and
\[
Q^N_n(\la,a):=(-\la)^{2N+1}+a_1\la^{2N}+...+a_{2N+1}\,, 
\]
with $a:=(a_1,...,a_{2N+1})$ in $\C^{2N+1}$.
For any $n\in\Z$, define $F^n:=(F^n_m)_{m\ne n}$ where for any $m\in\Z$, $m\ne n$,
\begin{equation}\label{e:F}
F^n_m(\si,a,\varphi):=(n-m)\oint_{\Ga_m}\frac{f_n(\la,\si,a)}{\sqrt[c]{\De(\la,\varphi)^2-4}}\,d\la\,.
\end{equation}
By construction, $F^n_m(\si,a,\varphi)$ is an analytic function on $l^2\times\C^{2N}\times\W$ for $|n|\le N$, and
an analytic function on $l^2\times\C^{2N+1}\times\W$ for $|n|\ge N+1$.
For $R,r>0$ consider the closed balls 
\[
B_R:=\{\si\in l^2\,|\,\|\si\|\le R\}\subseteq l^2
\]
and
\[
B_r^k:=\{a\in\C^k\,|\,|a|\le r\}\subseteq\C^k,
\]
where $k\ge 1$ is a given integer and $|a|=\sqrt{\sum_{j=1}^k|a_k|^2}$.
\begin{Lemma}\hspace{-2mm}{\bf .}\label{Lem:F-analytic}
For any $|n|\le N$ $[\mbox{resp.}\,\,|n|\ge N+1]$ the map
\[
l^2\times\C^{2N}\times\W\,\,[\mbox{resp.}\,\,\,l^2\times\C^{2N+1}\times\W]\to l^2,\,\,\,(\si,a,\varphi)\mapsto F^n(\si,a,\varphi)
\]
is well-defined and analytic. Moreover, for any $R,r>0$, and $|m|\ge N+1$
\begin{equation}\label{e:F-bounds}
F^n_m(\si,a,\varphi)=O(|\si_m|+|m\pi-\tau_m|+|\ga_m|)\,,
\end{equation}
uniformly in $B_R\times B_r^{2N}\times\W$ [resp. $B_R\times B_r^{2N+1}\times\W$]
and uniformly in $n\in\Z$ and $|m|\ge N+1$, $m\ne n$. In particular, by shrinking the neighborhood $\W$ if necessary,
$F^n$ is bounded in $B_R\times B_r^{2N}\times\W$ [resp. $B_R\times B_r^{2N+1}\times\W$] uniformly in $n\in\Z$.
\end{Lemma}
\noindent{\em Proof.}
First consider the case $|n|\ge N+1$.
For any $m\ne n$, the component $F^n_m$ of $F^n$ is analytic on $l^2\times\C^{2N+1}\times\W$ by construction.
Hence the analyticity of $F^n : l^2\times\C^{2N+1}\times\W\to l^2$ will follow once we prove that it is locally bounded
-- see e.g. Appendix A in \cite{KP}. For any $\si\in l^2$ and $|m|\ge N+1$ with $m\ne n$, one has
\begin{equation}\label{e:decomposition}
\frac{f_n(\la,\si,a)}{\sqrt[c]{\De(\la,\varphi)^2-4}}=
\frac{\tsi_m-\la}{\sqrt[s]{(\la-\la_m^-)(\la-\la_m^+)}}\,{\mathcal A}_n(\la,a,\varphi){\mathcal B}_m(\la,\si,\varphi)
\end{equation}
where, with $\tsi_n:=n\pi$,
\[
{\mathcal A}_n(\la,a,\varphi):=
i\,\frac{\la^{2N+1}+a_1\la^{2N}+...+a_{2N+1}}{(\tsi_n-\la)\prod\limits_{|j|\le N}\sqrt[s]{(\la-\la_j^-)(\la-\la_j^+)}}
\]
and
\[
{\mathcal B}_m(\la,\si,\varphi):=\prod\limits_{|j|\ge N+1,j\ne m}\frac{\tsi_j-\la}{\sqrt[s]{(\la-\la_j^-)(\la-\la_j^+)}}\,.
\]
A simple estimate shows that ${\mathcal A}_n(\la,a,\varphi)=O(1/|n-m|)$ uniformly on
$D_m\times B_r^{2N+1}\times\W$ and uniformly in $|m|\ge N+1$, $m\ne n$, and $|n|\ge N+1$.
By Lemma \ref{Lem:auxiliary3} in Appendix A, ${\mathcal B}_m(\la,\si,\varphi)=O(1)$ uniformly on
$D_m\times B_R\times\W$ and uniformly in $|m|\ge N+1$.
Combining these estimates for ${\mathcal A}_n$ and ${\mathcal B}_m$ with \eqref{e:decomposition} one gets from Lemma \ref{Lem:auxiliary4} in Appendix A that,
\[
F^n_m(\si,a,\varphi)=O(\rho_m)\,,\,\,\,\,\,\,\,\rho_m(\si):=\max\limits_{\la\in\Ga_m}|\tsi_m-\la|,
\]
uniformly on $B_R\times B_r^{2N+1}\times\W$ and uniformly in $|m|\ge N+1$, $m\ne n$, and $|n|\ge N+1$.
By shrinking the cycles $\Ga_m$ to $O(|\ga_m|)$-neighborhoods of $\tau_m$ one gets
from the triangle inequality that
\[
\rho_m(\si)=O(|\si_m|+|m\pi-\tau_m|+|\ga_m|)
\]
uniformly in $\si\in l^2$ and uniformly in $|m|\ge N+1$, $m\ne n$. This proves estimate \eqref{e:F-bounds}.
In a similar way one shows that for $|m|\le N$ and $|n|\ge N+1$, $F^n_m$ is bounded on $B_R\times B_r^{2N+1}\times\W$,
uniformly in $|m|\le N$ and $|n|\ge N+1$, implying that $F^n : B_R\times B_r^{2N+1}\times\W\to l^2$ is bounded,
uniformly in $|n|\ge N+1$. The case $|n|\le N$ is proved in a similar way.
\finishproof

\vspace{0.2cm}

\noindent\textit{Implicit function theorem:} 
Choose $N\ge N_0$, a compact set ${\mathcal K}\subseteq\C^{2N+1}$ as in Corollary \ref{Coro:K},
and $R>0$ so that 
\begin{equation}\label{e:enough_space1}
\|({\tau}_j(\varphi_*)-j\pi)_{|j|\ge N+1}\|<R/2\,.
\end{equation}
Furthermore, choose $r>0$ so that
\begin{equation}\label{e:enough_space2}
{\mathcal K}\subseteq B^{2N+1}_{r/2}
\end{equation}
and, for any $|n|\le N$, the coefficient vector $(a^{(n)}_j)_{1\le j\le 2N}$ of the polynomial
$P_n^N(\la,\varphi_*)=(-\la)^{2N}+a_1^{(n)}\la^{2N-1}+...+a_{2N}^{(n)}$ in \eqref{e:P_n^N} is contained in $B^{2N}_{r/2}$.
By shrinking the neighborhood $\W$ of $\varphi_*$ in $\LLb$, if necessary, we see from Lemma \ref{Lem:F-analytic} that
for any $|n|\le N$ $[\,\mbox{resp.}\,|n|\ge N+1]$ the analytic map
\[
l^2\times\C^{2N}\times\W\,[\,\mbox{resp.}\,\,l^2\times\C^{2N+1}\times\W]\to l^2,\,\,\,
(\si,a,\varphi)\mapsto F^n(\si,a,\varphi)
\]
is bounded on $B_R\times B_r^{2N}\times\W$ [resp. $B_R\times B_r^{2N+1}\times\W$] uniformly in $n\in\Z$.
For any $|n|\le N$ $[\,\mbox{resp.}\,|n|\ge N+1]$, denote by $a^{(n)}_*$ the coefficient vector
$(a_j^{(n)})_{1\le j\le 2N}$ [resp. $(a_j^{(n)})_{1\le j\le 2N+1}$] of the polynomial $P_n^N(\la,\varphi_*)$ and
let $\si_*^{(n)}:=(\tau_j(\varphi_*)-j\pi)_{|j|\ge N+1,j\ne n}$.
Note that by the choice of $r>0$ and Corollary \ref{Coro:K}, $a_*^{(n)}\in B_r^{2N}$ [resp. $a_*^{(n)}\in B_r^{2N+1}$]
for any $n\in\Z$. The main result of Section\,\,\ref{sec:finite_gaps} states that for any $n\in\Z$,
\begin{equation}\label{e:inverse_theorem_set_up}
F^n(\si_*^{(n)},a_*^{(n)},\varphi_*)=0\,.
\end{equation}
We now want to apply the implicit function theorem to show that for any $|n|\le N\,[\,\mbox{resp.}\,|n|\ge N+1]$,
there exist an open neighborhood $\W_n$ of $\varphi_*$ and analytic functions 
\[
\si^{(n)} : \W_n\to l^2,\,\,\,\,\varphi\mapsto\si^{(n)}(\varphi)
\]
and
\[
a^{(n)} : \W_n\to\C^{2N}\,[\,\mbox{resp.}\,\C^{2N+1}]
\]
such that for any $\varphi\in\W_n$, $F^n(\si^{(n)}(\varphi),a^{(n)}(\varphi),\varphi)=0$.
The function
\begin{equation}\label{e:zeta_n}
\zeta_n(\la,\varphi):=f_n(\si^{(n)}(\varphi),a^{(n)}(\varphi),\la)
\end{equation}
then has the required properties (see Corollary \ref{Coro:implicit_function} below).
In Section\,\,\ref{sec:uniformity} we then prove that one can choose the neighborhood $\W_n$
to be independent of $n\in\Z$.

\section{Differential of $F^n$}\label{sec:implicit_function}
Choose $N_0$ and $N\ge N_0$ as in Corollary \ref{Coro:K} and for any $n\in\Z$ let $F^n$ be the map introduced in
the previous section. In order to be able to apply the implicit function theorem to the equation $F^n(\si,a,\varphi)=0$
we show that the differential of $F^n$ w.r. to $(\si,a)$ at the point $(\si_*^{(n)},a_*^{(n)},\varphi_*)$ is a linear isomorphism.
It is convenient to denote the pair $(\si,a)$ by a single sequence $v=(v_j)_{j\ne n}$ by inserting for $|n|\le N$ [resp. $|n|\ge N+1$]
the vector $(a_j)_{1\le j\le 2N}$ [resp. $(a_j)_{1\le j\le 2N+1}$] in the middle of the sequence
$(\si_j)_{|j|\ge N+1,j\ne n}$. More precisely for $|n|\ge N+1$ and $k\in\Z\setminus\{n\}$ define
\[
v_k:=\left\{
\begin{array}{cl}
a_{N+1+k},&|k|\le N,\\
\si_k,&|k|\ge N+1.
\end{array}
\right. 
\]
Similarly, for $|n|\le N$ and $k\in\Z\setminus\{n\}$, define
\[
v_k:=\left\{
\begin{array}{cl}
a_{N+1+k},&-N\le k<n,\\
a_{N+1+k-1},&n<k\le N,\\
\si_k,&|k|\ge N+1.
\end{array}
\right. 
\]
Let $v^{(n)}_*$ be the sequence corresponding to $(\si_*^{(n)},a_*^{(n)})$.
We now compute the partial derivative $\partial_v F^n(v_*^{(n)},\varphi_*)$.
Recall that for any $m,n\in\Z$ with $m\ne n$
\begin{equation}
F^n_m(v,\varphi):=(n-m)\oint_{\Ga_m}\frac{f_n(\la,v)}{\sqrt[c]{\De(\la,\varphi)^2-4}}\,d\la\,.
\end{equation}
First consider the case $|n|\ge N+1$. Then by \eqref{e:f_n2},
$\frac{(n-m)f_n(\la,v)}{\sqrt[c]{\De(\la,\varphi)^2-4}}$ equals,
\[
\begin{array}{l}
\frac{i(n-m)}{\sqrt[s]{(\la-\la_n^-)(\la-\la_n^+)}}
\frac{Q_n^N(\la,v)}{\prod\limits_{|k|\le N}\sqrt[s]{(\la-\la_k^-)(\la-\la_k^+)}}
\prod\limits_{|k|\ge N+1,k\ne n}\frac{{\tilde\si}_k-\la}{\sqrt[s]{(\la-\la_k^-)(\la-\la_k^+)}}\,,
\end{array}
\]
where for any $v=(\si,a)$ one sets $Q_n^N(\la,v):=Q_n^N(\la,a)$.
Hence, by Cauchy's theorem, for any $|j|\ge N+1$ and $|m|\ge N+1$ with $m,j\ne n$
\begin{equation}\label{e:A)diagonal}
\partial_{v_j} F^n_m(v_*^{(n)},\varphi_*)=
\frac{2\,\pi(n-m)\,Q^N_n(\tau_{j},a_*^{(n)})}
{(\tau_{n}-\tau_{j})\prod\limits_{|k|\le N}\sqrt[s]{(\tau_{j}-\la_k^-)(\tau_{j}-\la_k^+)}}\,\de_{jm}\,,
\end{equation}
where $\tau_{k}\equiv\tau_k(\varphi_*)$ and $\la_k^\pm\equiv\la_k^\pm(\varphi_*)$ for any $k\in\Z$.
As 
\[
Q_n^N(\la,v)=(-\la)^{2N+1}+v_{-N}\la^{2N}+...+v_j\la^{N-j}+...+v_N
\]
one has for $|j|\le N$ and $m\ne n$,
\begin{equation}\label{e:A)central}
\partial_{v_j} F^n_m(v_*^{(n)},\varphi_*)=
(n-m)\,\oint_{\Ga_m}\frac{i\la^{N-j}}
{(\tau_{n}-\la)\prod\limits_{|k|\le N}\sqrt[s]{(\la-\la_k^-)(\la-\la_k^+)}}\,d\la\,.
\end{equation}
In particular, we get by Cauchy's theorem,
\begin{equation}\label{e:zeroes}
\partial_{v_j} F^n_m(v_*^{(n)},\varphi_*)=0\,\,\,\,\,\,\,\forall\,\,|j|\le N,\,\,|m|\ge N+1,\,\,m\ne n\,.
\end{equation}
Hence, \eqref{e:A)diagonal} holds for any $j\ne n$ and any $|m|\ge N+1$, $m\ne n$.
Later we will need also the case $|m|\le N$ and $|j|\ge N+1$ with $j\ne n$,
\begin{equation}\label{e:A)tail}
\partial_{v_j} F^n_m(v_*^{(n)},\varphi_*)=
\,\oint_{\Ga_m}\frac{i\,(n-m)\,Q_n^N(\la,a_*^{(n)})}
{(\tau_n-\la)(\tau_{j}-\la)\prod\limits_{|k|\le N}\sqrt[s]{(\la-\la_k^-)(\la-\la_k^+)}}\,d\la\,.
\end{equation}
Now, consider the case $|n|\le N$. Then by \eqref{e:f_n1}, $\frac{(n-m)f_n(\la,v)}{\sqrt[c]{\De(\la,\varphi)^2-4}}$ equals
\[
\begin{array}{l}
i(n-m)\frac{Q_n^N(\la,v)}{\prod\limits_{|k|\le N}\sqrt[s]{(\la-\la_k^-)(\la-\la_k^+)}}
\prod\limits_{|k|\ge N+1}\frac{{\tilde\si}_k-\la}{\sqrt[s]{(\la-\la_k^-)(\la-\la_k^+)}}
\end{array}
\]
and $Q_n^N(\la,v)=(-\la)^{2N}+v_{-N}\la^{2N-1}+...+v_{n-1}\la^{N-n}+v_{n+1}\la^{N-n-1}+...+v_N$.
Hence, for any $|j|\ge N+1$, $|m|\ge N+1$,
\begin{equation}\label{e:B)diagonal}
\partial_{v_j} F^n_m(v_*^{(n)},\varphi_*)=
\frac{2\,\pi(n-m)\,Q^N_n(\tau_{j},a_*^{(n)})}
{\prod\limits_{|k|\le N}\sqrt[s]{(\tau_{j}-\la_k^-)(\tau_{j}-\la_k^+)}}\,\de_{jm}\,.
\end{equation}
For $|j|\le N$ and $m,j\ne n$,
\begin{equation}\label{e:B)central}
\partial_{v_j} F^n_m(v_*^{(n)},\varphi_*)=
(n-m)\,\oint_{\Ga_m}\frac{i\,\la^{N-\de_j-1}}
{\prod\limits_{|k|\le N}\sqrt[s]{(\la-\la_k^-)(\la-\la_k^+)}}\,d\la
\end{equation}
where 
\[
\de_j:=\left\{
\begin{array}{cl}
j,&j<n,\\
j-1,&j>n.
\end{array}
\right.
\]
In particular, $\partial_{v_j} F^n_m(v_*^{(n)},\varphi_*)=0$
for $|j|\le N$, $|m|\ge N+1$. Hence \eqref{e:B)diagonal} holds for any $j\ne n$ and
$|m|\ge N+1$ with $m\ne n$.

\begin{Lemma}\hspace{-2mm}{\bf .}\label{Lem:dF-injective}
Let $\varphi_*\in\LLb$ be a finite gap potential. Then for any $n\in\Z$,
the linear map $\partial_v F^n(v_*^{(n)},\varphi_*)\in{\mathcal L}(l^2,l^2)$ is injective.
\end{Lemma}
{\em Proof.} Let us first treat the case $|n|\ge N+1$.
To simplify notation write $\partial_{v_j}F^n_m$ for $\partial_{v_j} F^n_m(v_*^{(n)},\varphi_*)$.
According to \eqref{e:A)diagonal} and \eqref{e:zeroes}, the infinite matrix $(\partial_{v_j}F^n_m)_{j,m\ne n}$ is in block form.
By \eqref{e:zeroes} the block $(\partial_{v_j}F^n_m)_{|m|, |j|\ge N+1, m,j\ne n}$ is diagonal and, by
Corollary \ref{Coro:K} $(i)$, none of the diagonal elements $\partial_{v_m}F^n_m$, $|m|\ge N+1, m\ne n$ vanishes. 
Further, by \eqref{e:zeroes}, the block $(\partial_{v_j}F^n_m)_{|m|\ge N+1, m\ne n, |j|\le N}$ vanishes.
Hence, to show that $\partial_v F^n : l^2\to l^2$ is one-to-one it suffices to show that the central block 
\[
M^n:=\Big(\partial_{v_j}F^n_m\Big)_{|m|\le N,|j|\le N}
\]
is non-degenerate. To prove it we want to apply Lemma \ref{Lem:a-periods} in Appendix B.
In the set-up of Appendix B choose $\Sigma$ be the compact Riemann surface $\Sigma_J$ of genus $g=|J|-1$ defined 
by \eqref{e:Sigma_J} in Section \ref{sec:finite_gaps},
with $J\equiv J(\varphi_*)\subseteq\Z$ and $\la_k^\pm=\la_k^\pm(\varphi_*)$, $k\in\Z$.
For $\Sigma^+$ we choose its canonical branch $\Sigma_J^c$, defined by \eqref{e:Sigma_J^c}.
The set of $C^1$-smooth simple closed curves $C_k$, $0\le k\le g$, in $(D1)$ of Appendix B is chosen to be
the set of curves $\pi_J^{-1}(G_j)$, $j\in J$, oriented according to $(D1)$.
Here $\pi_J : \CC_J\to\C$, denotes the projection $(\la,w)\mapsto\la$.
Next, introduce the holomorphic involution
\[
\imath : \CC_J\to\CC_J,\,\,\,\,\,\,(\la,w)\mapsto(\la,-w)
\]
and extend it to all of $\Sigma_J$ by setting $\imath(\infty^+)=\infty^-$ and $\imath(\infty^-)=\infty^+$.
As sets, $\imath(C_k)=C_k$ for any $0\le k\le g$. Finally, choose the points $\{P_k^+\}_{k=1}^s\subseteq\Sigma^+$, $s=2N+2-|J|$, so that
$\{P_k^+\}_{k=1}^{s-1}=\{\tau_{j}^+\,|\,j\notin J,|j|\le N\}$ and $P_s^+=\tau_{n}^+$ where $\pi_J(\tau_{j}^+)=\tau_{j}$.
The cycles $C_{g+k}$, $k=1,...,s$, are chosen as in $(D4)$ of Appendix B.
In view of formula \eqref{e:A)central}, for $1\le l\le 2N+1$, we define the meromorphic differentials on $\Sigma$
\begin{equation}\label{e:eta_l}
\eta_l=\frac{i\,\la^{2N-l+1}}{(\tau_{n}-\la)\prod\limits_{|k|\le N, k\notin J}(\tau_{k}-\la)\sqrt{\prod_{k\in J}(\la-\la_k^-)(\la-\la_k^+)}}\,d\la\,.
\end{equation}
In a straightforward way one checks that, for any $1\le l\le 2N+1$, $\eta_l$ is holomorphic at $\infty^{\pm}$ and may have poles
only at the points $P^\pm_k$, $1\le k\le s$, so that the divisor $(\eta_l)$ of $\eta_l$ satisfies
\[
(\eta_l)\ge-\sum\limits_{k=1}^s(P_k^-+P_k^+).
\]
Clearly, the differentials $\eta_l$, $1\le l\le 2N+1$, are linearly independent and satisfy $\imath(\eta_l)=-\eta_l$.
Hence, conditions $(D1)-(D5)$ of Appendix B are satisfied and one can apply Lemma \ref{Lem:a-periods}
with $k_0=g+s$ to conclude that the central block $M^n$ is non-degenerate.

The case $|n|\le N$ is treated similarly. The non-degeneracy of the corresponding central block
\[
M^n:=\Big(\partial_{v_j}F^n_m(v_*^{(n)},\varphi_*)\Big)_{|m|\le N,|j|\le N,\,m,j\ne n}
\]
follows from \eqref{e:B)central} and Lemma \ref{Lem:a-periods} with $k_0=g+s$, where $\Sigma$, $\Sigma^+$ and
the cycles $C_k$, $0\le k\le g$, are as above, but the points $\{P_k^+\}_{k=1}^s\subseteq\Sigma^+$, with $s=2N+1-|J|$,
are now chosen so that  $\{P_k^+\}_{k=1}^{s-1}=\{\tau_{j}^+\,|\,|j|\le N,j\notin J\cup\{n\}\}$ and
$P_s^+=\tau_{n}^+$. Finally, the cycles $C_{g+k}$, $1\le k\le s$, are chosen as in $(D4)$.
\finishproof

\begin{Prop}\hspace{-2mm}{\bf .}\label{Prop:dF-isomorphism}
Let $\varphi_*\in\LLb$ be a finite gap potential.
For any $n\in\Z$, $\partial_v F^n(v_*^{(n)},\varphi_*)\in{\mathcal L}(l^2,l^2)$
is an isomorphism.
\end{Prop}
{\em Proof.}
Consider the operator,
\[
D^n:=\left\{
\begin{array}{cl}
\partial_{v_j} F^n_m\,,&|m|\ge N+1\,\,\,\mbox{and}\,\,\, |j|\ge N+1\\
2\,\de_{jm},&|m|\le N\,\,\,\mbox{or}\,\,\,|j|\le N.
\end{array}
\right.
\]
In view of \eqref{e:A)diagonal} and \eqref{e:B)diagonal}, $D^n$ is a diagonal operator and by the choice of $N\ge N_0$
(so that Corollary \ref{Coro:K} holds) none of the diagonal entries of $D^n$ vanishes.
Using that by Proposition \ref{Prop:periodic_spectrum}, $\tau_{j}=j\pi+o(1)$ as $j\to\pm\infty$ one
gets from \eqref{e:A)diagonal} and \eqref{e:B)diagonal} that 
\begin{equation}\label{e:diagonal_limit}
\lim\limits_{m\to\pm\infty}D^n_{mm}=2\,.
\end{equation}
This implies that $D^n,(D^n)^{-1}\in{\mathcal L}(l^2,l^2)$.
Further, it follows from \eqref{e:A)diagonal}-\eqref{e:B)central} that the range of the operator
$K:=\partial_v F^n-D^n$ is contained in a finite dimensional space and hence $K$ is compact.
Note that,
\begin{equation}\label{e:fredholm}
\partial_v F^n=D^n\Big(\id+(D^n)^{-1}K\Big)
\end{equation}
where $\id$ is the identity on $l^2$.
As $\partial_v F^n$ and $D^n$ are injective we conclude from \eqref{e:fredholm} that $\id+(D^n)^{-1}K$ is injective. 
Therefore, by the Fredholm alternative, $\id+(D^n)^{-1}K$ is an isomorphism.
Combined with \eqref{e:fredholm}, it then follows that $\partial_v F^n$ is an isomorphism.
\finishproof

Proposition \ref{Prop:dF-isomorphism} allows to apply the implicit function theorem leading to
\begin{Coro}\hspace{-2mm}{\bf .}\label{Coro:implicit_function}
Let $\varphi_*\in\LLb$ be a finite gap potential. Then for any $n\in\Z$, there exist an open neighborhood $\W_n$ of $\varphi_*$
in $\LLb$ and an analytic function $\zeta_n : \C\times\W_n\to\C$ such that for any $\varphi\in\W_n$ and $m\in\Z$
\begin{equation}\label{e:normalization_local}
\frac{1}{2\pi}\oint_{\Ga_m}\frac{\zeta_n(\la,\varphi)}{\sqrt[c]{\De(\la,\varphi)^2-4}}\,d\la=\de_{mn}\,.
\end{equation}
In addition, for $\varphi=\varphi_*$, $\zeta_n(\la,\varphi_*)$ coincides with the entire function constructed in Theorem \ref{Th:main_finite_gap}.
\end{Coro}
\noindent{\em Proof.}
Let $|n|\le N$ $[\mbox{\rm resp.}\,|n|\ge N+1]$.
By Proposition \ref{Prop:dF-isomorphism} one can apply the implicit function theorem to
the analytic function  $F^n : l^2\times\C^{2N}\times\W\,[\,\mbox{resp.}\,l^2\times\C^{2N+1}\times\W]\to l^2$ constructed
in Section \ref{sec:set-up} to conclude that there exist an open neighborhood $\W_n$ of $\varphi_*$ in $\W$ and
analytic functions
\[
\si^{(n)} : \W_n\to l^2,\,\,\,\,\varphi\mapsto\si^{(n)}(\varphi)
\]
and
\[
a^{(n)} : \W_n\to\C^{2N}\,[\,\mbox{resp.}\,\C^{2N+1}]
\]
with $\si^{(n)}(\varphi_*)=\si^{(n)}_*$ and $a^{(n)}(\varphi_*)=a^{(n)}_*$ such that
for any $\varphi\in\W_n$, 
\[
F^n(\si^{(n)}(\varphi),a^{(n)}(\varphi),\varphi)=0\,.
\]
Hence, for any $\varphi\in\W_n$,
\[
\zeta_n(\la,\varphi):=f_n(\si^{(n)},a^{(n)},\varphi)
\]
satisfies
\[
\frac{1}{2\pi}\oint_{\Ga_m}\frac{\zeta_n(\la,\varphi)}{\sqrt[c]{\De(\la,\varphi)^2-4}}\,d\la=0\,\,\,\,\forall m\in\Z\setminus\{n\}\,.
\]
To see that the normalization condition
\begin{equation}\label{e:m=n}
\frac{1}{2\pi}\oint_{\Ga_n}\frac{\zeta_n(\la,\varphi)}{\sqrt[c]{\De(\la,\varphi)^2-4}}\,d\la=1
\end{equation}
holds we argue as follows. First note that by construction, $\zeta_n(\la,\varphi_*)$ coincides with the entire function
$\zeta_n$ of Theorem \ref{Th:main_finite_gap}. Thus in particular, \eqref{e:m=n} holds for $\varphi=\varphi_*$.
Furthermore, by \cite{KST2}, the finite gap potentials are dense in $\LL$. As $\LLb$ is open in $\LL$ the set of finite gap potentials
in $\LLb$ is also dense in $\LLb$. By continuity, it then suffices to prove \eqref{e:m=n} for finite gap potentials in $\W_n$.
It turns out that similar arguments as in the proof of Theorem \ref{Th:main_finite_gap} lead to the claimed result for finite gap potentials in $\W_n$. 
Indeed, for an arbitrary finite gap potential $\varphi\in\W_n$, denote by $J\equiv J(\varphi)\subseteq\Z$ the set of all $k\in\Z$ so that
$\la_k^-\ne\la_k^+$. First, consider the case when $n\notin J$. Then, by construction,
$\frac{\zeta_n(\la,\varphi)}{\sqrt[c]{\De(\la,\varphi)^2-4}}$ equals
\begin{equation}\label{e:decomposition2}
\begin{array}{l}
\frac{i}{\tau_n-\la}
\frac{Q_n^N(\la,a^{(n)})}
{\prod\limits_{|k|\le N,k\notin J,k\ne n}\!\!\!\!\!\!\!\!\!\!(\tau_k-\la)
\cdot\prod\limits_{k\in J}\sqrt[s]{(\la-\la_k^-)(\la-\la_k^+)}}
\prod\limits_{|k|\ge N+1,k\ne n}\!\!\!\!\frac{{\tilde\si}_k^{(n)}-\la}{\tau_k-\la}
\end{array}
\end{equation}
where ${\tilde\si}_k^{(n)}=k\pi+\si_k^{(n)}$.
As $\oint_{\Ga_m}\frac{\zeta_n(\la,\varphi)}{\sqrt[c]{\De(\la,\varphi)^2-4}}\,d\la=0$ for any $m\in\Z\setminus\{n\}$
one concludes that the the residue of \eqref{e:decomposition2} at $\tau_k$ vanishes for any $k\in\Z\setminus(J\cup\{n\})$.
Using that ${\tilde\si}_k^{(n)}=k\pi+o(1)$ as $k\to\pm\infty$ we get that
\[
\frac{\zeta_n(\la,\varphi)}{\sqrt[c]{\De(\la,\varphi)^2-4}}=
\frac{i}{\tau_n-\la}\frac{Q(\la)}{\prod\limits_{k\in J}\sqrt[s]{(\la-\la_k^-)(\la-\la_k^+)}}
\]
where $Q(\la)=(-\la)^{|J|}+\cdots$ is a polynomial of degree $|J|$.
Consider the Abelian differential
\[
\xi:=\frac{i}{\tau_n-\la}\frac{Q(\la)}{\sqrt{\prod\limits_{k\in J}(\la-\la_k^-)(\la-\la_k^+)}}\,d\la\,.
\]
It meromorphically extends to the compact Riemann surface $\Sigma_J$ introduced in Section \ref{sec:finite_gaps}
and has precisely four simple poles.
They are located at $\tau_n^{\pm}$ and $\infty^{\pm}$. Recall that $\pi_J(\tau_n^\pm)=\tau_n$ and
$\tau_n^+$ lies on the canonical branch $\Sigma_J^c$ of $\Sigma_J$ (see Section \ref{sec:finite_gaps}).
A straightforward computation shows that $\res\limits_{\infty^{\pm}}\xi=\pm i$. Finally, using that 
$\frac{1}{2\pi}\oint_{A_m}\xi=0$ for any $m\in J$ and that by Stokes' formula
\[
\sum_{m\in J}\oint_{A_m}\xi+2\pi i\res_{\tau_n^+}\xi+2\pi i\res_{\infty^+}\xi=0
\]
one concludes that $\res\limits_{\tau^+_n}\xi=-i$.
This completes the proof of the normalization condition \eqref{e:m=n} when $n\notin J$.
The case $n\in J$ is treated similarly.
\finishproof

\section{Uniformity}\label{sec:uniformity}
To prove Theorem \ref{Th:main} it remains to show that the neighborhoods $\W_n$ of $\varphi_*\in\LLb$ of
Corollary \ref{Coro:implicit_function} can be chosen independently of $n\in\Z$. 
Choose $N_0$ and $N\ge N_0$ as in the implicit function theorem part of Section \ref{sec:set-up}. 
First note that it suffices to show that $\W_n$ can be chosen independently of $n$ for any $|n|\ge N+1$.
For this purpose we want to study the asymptotics
of $(\partial_{v_j}F^n_m)_{m,j}$ for $|n|$ large. To this end introduce for $a=(a_j)_{1\le j\le 2N+1}$ in $\C^{2N+1}$
the infinite matrix
\[
\F^\infty(a):=(\F^\infty_{mj}(a))_{m,j\in\Z}
\]
where, for $|m|\ge N+1$ and $j\in\Z$,
\begin{equation}\label{e:A)diagonal*}
\F^\infty_{mj}(a):=
\frac{2\,Q^N_n(\tau_{j},a)}
{\prod\limits_{|k|\le N}\sqrt[s]{(\tau_{j}-\la_k^-)(\tau_{j}-\la_k^+)}}\,\de_{jm}\,,
\end{equation}
for $|m|\le N$ and $|j|\le N$,
\begin{equation}\label{e:A)central*}
\F^\infty_{mj}(a):=\frac{i}{\pi}
\oint_{\Ga_m}\frac{\la^{N-j}}
{\prod\limits_{|k|\le N}\sqrt[s]{(\la-\la_k^-)(\la-\la_k^+)}}\,d\la\,,
\end{equation}
and, for $|m|\le N$ and $|j|\ge N+1$,
\begin{equation}\label{e:A)tail*}
\F^\infty_{mj}(a):=\frac{i}{\pi}
\oint_{\Ga_m}\frac{Q^N_n(\la,a)}
{(\tau_{j}-\la)\prod\limits_{|k|\le N}\sqrt[s]{(\la-\la_k^-)(\la-\la_k^+)}}\,d\la\,.
\end{equation}
Here $Q^N_n(\la,a)=(-\la)^{2N+1}+a_1\la^{2N}+...+a_{2N+1}$, $\la_k^\pm\equiv\la_k^\pm(\varphi_*)$,
and $\tau_{k}\equiv\tau_k(\varphi_ *)$.

\vspace{0.2cm}

We remark that the entries $\F^\infty_{mj}(a)$ with $m,j\ne n$ in \eqref{e:A)diagonal*}, \eqref{e:A)central*},
and \eqref{e:A)tail*} are formally obtained from \eqref{e:A)diagonal}, \eqref{e:A)central},
and \eqref{e:A)tail} by replacing $a_*^{(n)}$ with $a$ and then taking the limit as $n\to\infty$.

Let ${\mathcal K}\subseteq\C^{2N+1}$ be the compact set introduced in Corollary \ref{Coro:K}.
\begin{Lemma}\hspace{-2mm}{\bf .}\label{Lem:F^*-isomorphism}
For any $a\in{\mathcal K}$, the operator $\F^\infty(a)\in{\mathcal L}(l^2,l^2)$ is an isomorphism.
Moreover, there exists $0<C<\infty$ so that $\|\F^\infty(a)^{-1}\|\le C$ for any $a\in{\mathcal K}$.
\end{Lemma}
{\em Proof.}
Arguing as in the proof of Lemma \ref{Lem:dF-injective} and taking into account formulas 
\eqref{e:A)diagonal*} and \eqref{e:A)central*} one concludes from Corollary \ref{Coro:K} and Lemma \ref{Lem:a-periods}
that for any $a\in{\mathcal K}$, the map $\F^\infty(a) : l^2\to l^2$ is injective.
More precisely we apply Lemma \ref{Lem:a-periods} with $\Sigma^+$ being
the canonical branch of the compact Riemann surface $\Sigma:=\Sigma_J$, introduced in Section \ref{sec:finite_gaps}, and
$\{P_k^\pm\}_{k=1}^s$ being points in $\Sigma$ with $s=2N+2-|J|$, so that 
$\{P_k^\pm\}_{k=1}^{s-1}=\{\tau_j^\pm(\varphi_*)\,|\,j\notin J,|j|\le N\}$ and
$P_s^\pm=\infty^\pm$.
Arguing as in the proof of Proposition \ref{Prop:dF-isomorphism} one concludes that $\F^\infty(a)$ is
a linear isomorphism for any $a\in{\mathcal K}$.
As 
\[
\C^{2N+1}\to{\mathcal L}(l^2,l^2),\,\,\,\,\,\,a\mapsto\F^\infty(a)\,,
\]
is continuous\footnote{This follows directly from \eqref{e:A)diagonal*}-\eqref{e:A)tail*}.} and ${\mathcal K}\subseteq\C^{2N+1}$ is compact, it follows that there exists $C>0$ so that
$\|\F^\infty(a)^{-1}\|\le C$ for any $a\in{\mathcal K}$.
\finishproof

In the sequel we also need to consider certain restrictions of the operator $\F^\infty(a)$.
For any $|n|\ge N+1$, denote by $\F^\infty_{,n}(a)$ the restriction of $\F^\infty(a)$ to $l^2(\Z\setminus\{n\},\C)$,
\[
\F^\infty_{,n}(a) : l^2(\Z\setminus\{n\},\C)\to l^2(\Z\setminus\{n\},\C),\,\,\,\,\,\,\,
(\xi_l)_{l\ne n}\mapsto\Big(\sum\limits_{l\ne n}\F^\infty_{ml}\xi_l\Big)_{m\ne n}\,.
\]
Using the block structure of $\F^\infty(a)$ one easily gets
\begin{Coro}\hspace{-2mm}{\bf .}\label{Coro:F_{,n}}
For any $|n|\ge N+1$ and any $a\in{\mathcal K}$, $\F^\infty_{,n}(a)$ is a linear isomorphism and 
\[
(\F^\infty_{,n}(a))^{-1}=({\F^\infty(a)}^{-1})_{,n}\,.
\]
\end{Coro}
Combining Lemma \ref{Lem:F^*-isomorphism} with Corollary \ref{Coro:F_{,n}} we get
\begin{Coro}\hspace{-2mm}{\bf .}\label{Coro:F^*_{,n}-isomorphism}
For any $a\in{\mathcal K}$ and any $|n|\ge N+1$, the operator $\F^\infty_{,n}(a)\in{\mathcal L}(l^2,l^2)$ is
an isomorphism. Furthermore, there exists $0<C<\infty$ so that for any $a\in{\mathcal K}$ and $|n|\ge N+1$,
\[
\|\F^\infty_{,n}(a)^{-1}\|\le C\,.
\]
\end{Coro}
For any $|n|\ge N+1$, define
\begin{equation}\label{e:F_a}
\F^n : \C^{2N+1}\to{\mathcal L}(l^2,l^2),\,\,\,\,\,\,\,
a\mapsto\Big(\partial_{v_j}F^n_m(\si^{(n)}_*,a,\varphi_*)\Big)_{j,m\in\Z\setminus\{n\}}
\end{equation}
where $\si_*^{(n)}=(\tau_k-k\pi)_{|k|\ge N+1,k\ne n}$.
\begin{Lemma}\hspace{-2mm}{\bf .}\label{Lem:F^n->F^*}
For $a\in{\mathcal K}$, 
\[
\F^n(a)-\F^\infty_{,n}(a)\to 0\,\,\,\,\,\mbox{as}\,\,\,\ n\to\pm\infty
\]
in ${\mathcal L}(l^2,l^2)$, uniformly on ${\mathcal K}$.
\end{Lemma}
{\em Proof.} Assume that $|n|\ge N+1$. Arguing as for the derivation of
\eqref{e:A)diagonal}, \eqref{e:A)central}, and \eqref{e:A)tail} one gets
for any $|m|\ge N+1$, $m\ne n$, and $j\in\Z\setminus\{n\}$
\begin{equation}\label{e:A)diagonal_a}
\F^n_{mj}(a)=
\frac{2\,\pi(n-m)\,Q^N_n(\tau_{j},a)}
{(\tau_{n}-\tau_{j})\prod\limits_{|k|\le N}\sqrt[s]{(\tau_{j}-\la_k^-)(\tau_{j}-\la_k^+)}}\,\de_{jm}
\end{equation}
where $Q^N_n(\la,a)=(-\la)^{2N+1}+a_1\la^{2N}+...+a_{2N+1}$;
for $|m|\le N$ and $|j|\le N$,
\begin{equation}\label{e:A)central_a}
\F^n_{mj}(a)=
(n-m)\,\oint_{\Ga_m}\frac{i\,\la^{N-j}}
{(\tau_{n}-\la)\prod\limits_{|k|\le N}\sqrt[s]{(\la-\la_k^-)(\la-\la_k^+)}}\,d\la\,,
\end{equation}
whereas for $|m|\le N$ and $|j|\ge N+1$ with $j\ne n$,
\begin{equation}\label{e:A)tail_a}
\F^n_{mj}(a)=
(n-m)\,\oint_{\Ga_m}\frac{i\,Q^N_n(\la,a)}
{(\tau_{n}-\la)(\tau_{j}-\la)\prod\limits_{|k|\le N}\sqrt[s]{(\la-\la_k^-)(\la-\la_k^+)}}\,d\la\,.
\end{equation}
Using \eqref{e:A)tail*}, \eqref{e:A)tail_a}, and the asymptotic formula for $\tau_n$,
$\tau_n=n\pi+o(1)$, $n\to\pm\infty$, one gets for $|m|\le N$, $|j|\ge N+1$ with $j\ne n$, and $|n|\ge N+1$
\begin{eqnarray}
|\F^n_{mj}-\F^\infty_{mj}|&=&\frac{1}{\pi}
\left|\oint_{\Ga_m}\frac{\Big(\frac{\pi(n-m)}{\tau_{n}-\la}-1\Big)Q^N_n(\la,a)}
{(\tau_{j}-\la)\prod\limits_{|k|\le N}\sqrt[s]{(\la-\la_k^-)(\la-\la_k^+)}}\,d\la\right|\nonumber\\
&=&O\Big(\frac{1}{(n-m)(j-m)}\Big)\label{e:ineq1}\,,
\end{eqnarray}
uniformly in $|n|\ge N+1$, $|m|\le N$, $|j|\ge N+1$, and $a\in{\mathcal K}$.
Similarly, using \eqref{e:A)central*} and \eqref{e:A)central_a} and the asymptotics of $\tau_n$ as $n\to\pm\infty$
one gets for $|m|\le N$, $|j|\le N$, and $|n|\ge N+1$
\begin{equation}\label{e:ineq3}
|\F^n_{mj}-\F^\infty_{mj}|=O\Big(\frac{1}{n-m}\Big)
\end{equation}
uniformly in $a\in{\mathcal K}$ and $|n|\ge N+1$.
Finally, \eqref{e:A)diagonal*}, \eqref{e:A)diagonal_a}, and the asymptotics of $\tau_n$ imply that uniformly in
$|n|\ge N+1$, $|m|\ge N+1$, $m\ne n$, and $a\in{\mathcal K}$,
\begin{eqnarray}
|\F^n_{mm}-\F^\infty_{mm}|=O\Big(\frac{r_n+r_m}{|n-m|}\Big),\label{e:ineq2}
\end{eqnarray}
where for $|j|\ge N+1$, $r_j:=|\tau_j-j\pi|$.
Furthermore, for any $|m|\le N$, one gets from \eqref{e:ineq1} and \eqref{e:ineq3} that 
\begin{equation}\label{e:ineqA}
\|\F^n_{m\bullet}-\F^\infty_{m\bullet}\|=
\Big(\sum_{j\in\Z,j\ne n}|\F^n_{mj}-\F^\infty_{mj}|^2\Big)^{1/2}=O\Big(\frac{1}{|n-N|}\Big)\,
\end{equation}
uniformly in $|n|\ge N+1$ and $a\in{\mathcal K}$.
Using \eqref{e:ineq2}, one concludes that
\[
\sup_{|m|\ge N+1,m\ne n}|\F^n_{mm}-\F^\infty_{mm}|=O\Big(r_n+\max_{j\ne 0}\Big(\frac{r_{n-j}}{|j|}\Big)\Big)
\]
uniformly in $|n|\ge N+1$ and $a\in{\mathcal K}$.
Note that,
\begin{eqnarray*}
\max_{j\ne 0}\Big(\frac{r_{n-j}}{|j|}\Big)&\le&\max_{0<|j|\le n/2}\Big(\frac{r_{n-j}}{|j|}\Big)+
\max_{|j|> n/2}\Big(\frac{r_{n-j}}{|j|}\Big)\\
&\le&\max_{|j|\ge n/2}r_j+2\max_{j\in\Z}r_j/n=o(1)
\end{eqnarray*}
as $n\to\pm\infty$.
Combining the estimates obtained it follows that
\begin{equation}\label{e:ineqB}
\sup_{|m|\ge N+1,m\ne n}|\F^n_{mm}-\F^\infty_{mm}|=o(1),\,\,\,\,\,\,n\to\pm\infty
\end{equation}
uniformly in $a\in{\mathcal K}$. The claimed estimate now follows from \eqref{e:ineqA} and \eqref{e:ineqB}.
\finishproof

As an immediate consequence of Lemma \ref{Lem:F^n->F^*} and Corollary \ref{Coro:F^*_{,n}-isomorphism} one obtains
\begin{Coro}\hspace{-2mm}{\bf .}\label{Coro:bounded_inverse}
There exist $0<C<\infty$ and $N_1\ge N$ so that for any $|n|\ge N_1$ and $a\in{\mathcal K}$,
$\F^n(a) : l^2\to l^2$  is a linear isomorphism and
\begin{equation}\label{e:bounded_inverse}
\|\F^n(a)^{-1}\|\le C
\end{equation}
uniformly on ${\mathcal K}$.
\end{Coro}
\noindent{\em Proof of Theorem \ref{Th:main}.}
Choosing $N$ greater, if necessary, we obtain from Corollary \ref{Coro:bounded_inverse} that
\eqref{e:bounded_inverse} holds for any $|n|\ge N+1$.
In view of Corollary \ref{Coro:implicit_function} it remains to be proved that the neighborhoods
$\W_n$ of $\varphi_*$ in $\W$ with $|n|\ge N+1$ can be chosen independently of $n$.

Recall that for any $|n|\ge N+1$ the sequence $\si^{(n)}_*=(\tau_j(\varphi_*)-j\pi)_{|j|\ge N+1,j\ne n}$
belongs to the space
\[
l^2_{N,n}:=\{x=(x_j)_{|j|\ge N+1,j\ne n}\,|\,\|x\|_{N,n}<\infty\}
\]
where $\|x\|_{N,n}:=\Big(\sum_{|j|\ge N+1,j\ne n}|x_j|^2\Big)^{1/2}$.
Let
\[
l^2_N:=\{x=(x_j)_{|j|\ge N+1}\,|\,\|x\|_N<\infty\}
\]
with $\|x\|_N:=\Big(\sum_{|j|\ge N+1}|x_j|^2\Big)^{1/2}$.
For any $n\le -N-1$ consider the linear isomorphism $\imath_n : l^2_{N,n}\to l^2_N$, defined for any
$x=(x_j)_{|j|\ge N+1,j\ne n}\in l^2_{N,n}$ by
\[
(\imath_n(x))_j=
\left\{
\begin{array}{cc}
x_{j-1},&j\le n\\
x_j,&j>n
\end{array}
\right.
\]
and similarly for $n\ge N+1$,
\[
(\imath_n(x))_j=
\left\{
\begin{array}{cc}
x_j,&j<n\\
x_{j+1},&j\ge n
\end{array}
\right. .
\]
Clearly, for any $|n|\ge N+1$, $\imath_n$ is an isometry, i.e., 
\[
\|\imath_n(x)\|_N=\|x\|_{N,n}\,.
\] 
Using the isometries $\imath_n$ we identify  $l^2_{N,n}$ with $l^2_N$ for any $|n|\ge N+1$ and simply write
$\sigma^{(n)}_*\in l^2_N$. The functions $F^n$ are then analytic and bounded uniformly in $|n|\ge N+1$ on the set 
\[
B_R\times B_r^{2N+1}\times\W\subseteq l^2_N\times\C^{2N+1}\times\LL
\]
where now $B_R=\{x\in l^2_N\,|\,\|x\|_N\le R\}$.
By construction -- see \eqref{e:enough_space1}, \eqref{e:enough_space2} -- for any $|n|\ge N+1$
\begin{equation}\label{e:enough_space}
\|\si^{(n)}_*\|_N\le R/2\,\,\,\,\,\mbox{and}\,\,\,\,\,{\mathcal K}\subseteq B_{r/2}^{2N+1}\,.
\end{equation}
Denote by $|\!|\!|\cdot|\!|\!|$ the canonical norm of the Cartesian product, $l^2_N\times\C^{2N+1}\times\LL$, and let 
\[
z_n:=(\si^{(n)}_*,a^{(n)}_*,\varphi_*)
\]
for $|n|\ge N+1$.
By Corollary \ref{Coro:K} $(ii)$, $a^{(n)}_*\in{\mathcal K}$. Hence by \eqref{e:enough_space}, one can choose
$\rho>0$ such that for any $|n|\ge N+1$, the ball in $l^2_N\times\C^{2N+1}\times\LL$ of radius $\rho$, centered at $z_n$,
\[
B_\rho(z_n):=\{z\in l^2_N\times\C^{2N+1}\times\LL\,:\,\,|\!|\!|z-z_n|\!|\!|\le\rho\}\,,
\]
is contained in $B_R\times B_r^{2N+1}\times\W$.
As $F^n$ is bounded on $B_R\times B_r^{2N+1}\times\W$ uniformly in $|n|\ge N+1$, we obtain by Cauchy's estimate
(cf. \cite[Lemma A.2, Appendix A]{KP}) that the ${\mathcal L}(l^2,l^2)$-norm of
the derivative $\partial_v F^n$ is bounded on $B_{\rho/2}(z_n)$ by a constant independent of $|n|\ge N+1$.
Applying Cauchy's estimate once more we see that the analytic map
\[
B_R\times B_r^{2N+1}\times\W\to{\mathcal L}(l^2,l^2),\,\,\,(\si,a,\varphi)\mapsto (\partial_v F^n)(\si,a,\varphi)
\]
is Lipschitz continuous on $O_{\rho/4}(z_n)$ with a Lipschitz constant independent of $|n|\ge N+1$.
Together with Corollary \ref{Coro:bounded_inverse} one concludes that there exist $0<C<\infty$ and
$0<\rho_1<\rho/4$ such that for any $|n|\ge N+1$ and $z\in O_{\rho_1}(z_n)$,
\begin{equation}\label{e:uniform_estimates}
\|(\partial_vF^n(z))^{-1}\|\le 2C<\infty\quad\mbox{and}\quad\|\partial_v F^n(z)\|\le 2C<\infty\,.
\end{equation}
It then follows from the uniform estimates \eqref{e:uniform_estimates} and the implicit function theorem
that one can choose the neighborhoods $\W_n$ in Corollary \eqref{Coro:implicit_function} independently of $|n|\ge N+1$.
Denote this neighborhood by ${\widetilde\W}$ and set $\W:={\widetilde\W}\bigcap\Big(\cap_{|n|\le N}\W_n\Big)$.
By construction, $\zeta_n : \C\times\W\to\C$,
\begin{equation}\label{e:zeta_n'}
\zeta_n(\la,\varphi)=-\frac{2}{\pi_n}\prod_{k\ne n}\frac{\tsi_k^{(n)}-\la}{\pi_k}\,,\,\,\,\,\,\,\,
\tsi_k^{(n)}=\tsi_k^{(n)}(\varphi),
\end{equation}
and 
\begin{equation}\label{e:the_equation}
\frac{1}{2\pi}\oint_{A_m}\frac{\zeta_n(\la,\varphi)}{\sqrt{\De(\la,\varphi)^2-4}}\,d\la=\de_{mn}\,\,\,\,
\forall m\in\Z\,.
\end{equation}
In addition, for any $n\in\Z$, the map $\si^{(n)}: \W\to l^2$,
\begin{equation}\label{e:lipschitz}
\varphi\mapsto(\si^{(n)}_k(\varphi))_{|k|\ge N+1,k\ne n}\,\,\,\,\,\mbox{where}\,\,\,\,\,\si^{(n)}_k=\tsi_k^{(n)}-k\pi
\end{equation}
is analytic and the roots $\tsi_k^{(n)}$, $|k|\le N$, of $\zeta_n$ are contained in the disk
$\{\la\in\C\,|\,|\la|\le (N+\frac{1}{4})\pi\}$.

\vspace{0.2cm}

Now, we will prove the uniform estimate \eqref{e:asymptotics*}. 
For any $m\ne n$, equation \eqref{e:the_equation} can be written as
\begin{equation}\label{e:the_equation'}
\oint_{\Ga_m}\frac{\tsi^{(n)}_m-\la}{\sqrt[s]{(\la-\la_m^-)(\la-\la_m^+)}}\,\chi^n_m(\la,\varphi)\,d\la=0
\end{equation}
where
\begin{equation}\label{e:the_equation''}
\chi^n_m(\la,\varphi):=\frac{\pi(n-m)}{\tsi^{(n)}_n-\la}\,\prod_{j\ne m}\frac{\tsi^{(n)}_j-\la}{\sqrt[s]{(\la-\la_j^-)(\la-\la_j^+)}}
\end{equation}
with $\tsi^{(n)}_n:=n\pi$.
By shrinking the neighborhood $\W$, if necessary, we conclude from Lemma \ref{Lem:auxiliary3} that the function $\chi^n_m$
is bounded on $\D_m\times\W$ by a constant independent of $n\in\Z$ and $|m|\ge N+1$, $m\ne n$.
In addition, by shrinking the neighborhood $\W$ once more and choosing $N\ge N_0$ greater if necessary we can 
ensure, using Lemma \ref{Lem:some_analytic_maps}, that for any $|m|\ge N+1$ and $\varphi\in\W$,
\[
\tau_m,\,\la_m^\pm\in\{\la\in\C\,|\,|\la-m\pi|\le\pi/8\}\,.
\]
Assume that $|m|\ge N+1$, $m\ne n$. For any $(\la,\varphi)\in\D_m\times\W$ we have
\begin{equation}\label{e:taylor}
\chi^n_m(\la,\varphi)=\chi^n_m(\tau_m,\varphi)+(\la-\tau_m)g^n_m(\la,\varphi)
\end{equation}
where $g^n_m(\la):=\int_0^1{\dot\chi}^n_m(\tau_m+(\la-\tau_m)t)\,dt$ and ${\dot\chi}^n_m$ denotes
the derivative of $\chi^n_m$ with respect to $\la$. Using that $\chi^n_m$ is bounded on $\D_m\times\W$ uniformly in
$n\in\Z$ and $|m|\ge N+1$, $m\ne n$, Cauchy's estimate for ${\dot\chi}^n_m$ implies that $g^n_m$ is bounded on
$(\la,\varphi)\in\{\la\in\C\,|\,|\la-m\pi|\le\pi/8\}\times\W$ uniformly in $n\in\Z$ and $|m|\ge N+1$, $m\ne n$.
As 
\[
\frac{1}{2\pi i}\oint_{\Ga_m}\frac{\tsi^{(n)}_m-\la}{\sqrt[s]{(\la-\la_m^-)(\la-\la_m^+)}}\,d\la=\tau_m-\tsi^{(n)}_m
\]
formulas \eqref{e:the_equation'}-\eqref{e:taylor} then lead to
\begin{equation}\label{e:the_equation*}
(\tau_m-\tsi^{(n)}_m)\chi^n_m(\tau_m,\varphi)+\frac{1}{2\pi i}
\oint_{\Ga_m'}\frac{(\tsi^{(n)}_m-\la)(\la-\tau_m)}{\sqrt[s]{(\la-\la_m^-)(\la-\la_m^+)}}\,g^n_m(\la,\varphi)\,d\la=0
\end{equation}
where $\Ga_m':=\{\la\in\C\,|\,|\la-m\pi|=\pi/8\}$.
It follows from Lemma \ref{Lem:auxiliary3'} and Lemma \ref{Lem:some_analytic_maps} that by shrinking
the neighborhood $\W$, if necessary, once more we can find $c>0$ and $m_0>0$ such that
\begin{equation}\label{e:bounded_below}
|\chi_m^n(\tau_m,\varphi)|>c>0
\end{equation}
uniformly in $\varphi\in\W$, $n\in\Z$ and $|m|\ge N+1$, $m\ne n$. 
As by construction, $\|(\si_m^{(n)})_{|m|\ge N+1}\|<R$, Lemma \ref{Lem:auxiliary4} shows that the integral in
\eqref{e:the_equation*} is of order $O(\ga_m)$ uniformly for $\varphi\in\W$, $n\in\Z$, and $|m|\ge N+1$, $m\ne n$.
This together with \eqref{e:the_equation*} and \eqref{e:bounded_below} implies that $\tau_m-\tsi^{(n)}_m=O(\ga_m)$
uniformly in $\varphi\in\W$, $n\in\Z$, and $|m|\ge N+1$, $m\ne n$. 
Using Lemma \ref{Lem:auxiliary4} once more we then conclude that
uniformly for $\varphi\in\W$, $n\in\Z$ and $|m|\ge N+1$, $m\ne n$,
\begin{equation}\label{e:integral_estimate}
\oint_{\Ga_m'}\frac{(\tsi^{(n)}_m-\la)(\la-\tau_m)}{\sqrt[s]{(\la-\la_m^-)(\la-\la_m^+)}}\,g^n_m(\la,\varphi)\,d\la
=O(\ga_m^2)\,.
\end{equation}
Estimate \eqref{e:asymptotics*} then follows from \eqref{e:the_equation*}, \eqref{e:integral_estimate},
and \eqref{e:bounded_below}.

\noindent The last statement of the theorem follows directly from \eqref{e:canonical_root}, \eqref{e:zeta_n'},
\eqref{e:the_equation}, and Cauchy's formula.
\finishproof

\section{Appendix A: Estimates on products} 
In this Appendix we collect some technical lemmas used in the main body of the paper.
\begin{Lemma}\hspace{-2mm}{\bf .}\label{Lem:auxiliary1}
For $a=(a_j)_{j\in\Z}\in l^1$,
\begin{equation}\label{e:auxiliary}
\Big|\prod_{j\in\Z}(1+a_j)-1\Big|\le\|a\|_{l^1}\exp(\|a\|_{l^1})\,,
\end{equation}
where $\|a\|_{l^1}=\sum_{j\in\Z}|a_j|$.
\end{Lemma}
\noindent{\em Proof.}
As $a\in l^1$, the product $\prod\limits_{j\in\Z}(1+a_j)$ converges absolutely.
If $a_j=0$ for all $j\in\Z$ except finitely many then by the triangle inequality,
\begin{equation}\label{e:triangle}
\Big|\prod_{j\in\Z}(1+a_j)-1\Big|\le\prod_{j\in\Z}(1+|a_j|)-1\,.
\end{equation}
By a limiting argument one sees that this inequality holds also in the case of an arbitrary element $a\in l^1$.
As $\log(1+x)\le x$ for $x\ge 0$, one gets 
\begin{equation}\label{e:log}
\prod_{j\in\Z}(1+|a_j|)=\exp\Big(\sum_{j\in\Z}\log(1+|a_j|)\Big)\le\exp(\|a\|_{l^1})\,.
\end{equation}
By the Taylor expansion of $e^x$ at zero, one has
\[
0\le e^x-1\le x e^x\,\,\,\forall x\ge 0,
\]
implying \eqref{e:auxiliary}, by combining \eqref{e:triangle} and \eqref{e:log}.
\finishproof

\begin{Lemma}\hspace{-2mm}{\bf .}\label{Lem:auxiliary2}
For $\si=(\si_j)_{j\in\Z}\in l^2$ and any $m\in\Z$,
\[
\prod_{j\ne m}\Big(1+\frac{\si_j}{j-m}\Big)=1+r_m(\si)\,,
\]
where 
\begin{equation}\label{e:r_m-rough}
|r_m(\si)|\le 2\,\|\si\|\exp(2\|\si\|)
\end{equation}
and $\|\si\|=\Big(\sum_{j}|\si_j|^2\Big)^{1/2}$.
Moreover, for $|m|\ge 2$,
\begin{equation}\label{e:r_m}
|r_m(\si)|\le 2\,\Big(\frac{\|\si\|}{(|m|-1)^{1/2}}+\|T_m(\si)\|\Big)\exp(2\|\si\|)
\end{equation}
with $T_m(\si):=(\si_j)_{|j|\ge|m|/2}$.
\end{Lemma}
\noindent{\em Proof.}
Consider the sequence $a^{(m)}=(a^{(m)}_j)_{j\in\Z}$,
\[
a^{(m)}_j:=\left\{\begin{array}{cc}\frac{\si_j}{j-m}&j\ne m\\0&j=m\end{array}\right..
\]
By the Cauchy-Schwarz inequality,
\begin{equation}\label{e:a1}
\|a^{(m)}\|_{l^1}\le\|\si\|\Big(\sum_{j\ne m}\frac{1}{|j-m|^2}\Big)^{1/2}\le 2\|\si\|\,.
\end{equation}
Combining this with Lemma \ref{Lem:auxiliary1} we obtain \eqref{e:r_m-rough}.
In fact, \eqref{e:a1} can be improved,
\begin{eqnarray}
\|a^{(m)}\|_{l^1}&\le&\|\si\|\Big(\sum_{|j-m|>|m|/2}\frac{1}{|j-m|^2}\Big)^{1/2}\nonumber\\
&+&\Big(\sum_{1\le|j-m|\le|m|/2}\frac{1}{|j-m|^2}\Big)^{1/2}\Big(\sum_{|j|\ge |m|/2}|\si_j|^2\Big)^{1/2}\label{e:a2}\,.
\end{eqnarray}
For $|m|\ge 2$,
\begin{eqnarray}
\Big(\sum_{|j-m|>|m|/2}\frac{1}{|j-m|^2}\Big)^{1/2}&=
&\Big(2\sum_{k\ge\frac{|m|}{2}+\frac{1}{2}}\frac{1}{k^2}\Big)^{1/2}\nonumber\\
&\le&\Big(2\sum_{k\ge\frac{|m|}{2}+\frac{1}{2}}\frac{1}{(k-1)k}\Big)^{1/2}\nonumber\\
&\le&\Big(\frac{4}{|m|-1}\Big)^{1/2}\label{e:a3}
\end{eqnarray}
and
\begin{equation}\label{e:a4}
\Big(\sum_{1\le|j-m|\le|m|/2}\frac{1}{|j-m|^2}\Big)^{1/2}\le\Big(2\,\sum_{k\ge 1}\frac{1}{k^2}\Big)^{1/2}\le 2.
\end{equation}
Combining \eqref{e:a2} with \eqref{e:a3} and \eqref{e:a4} we get,
\begin{equation}\label{e:a5}
\|a^{(m)}\|_{l^1}\le 2\,\Big(\frac{\|\si\|}{(|m|-1)^{1/2}}+\|T_m(\si)\|\Big)\,.
\end{equation}
Finally, \eqref{e:r_m} follows from Lemma \ref{Lem:auxiliary1}, \eqref{e:a1}, and \eqref{e:a5}.
\finishproof

Let $\varphi_*\in\LLb$. Choose an open neighborhood $\W$ of $\varphi_*$ in $\LLb$, $N_0\ge 1$,
and cycles $\Ga_m$ ($m\in\Z$) as in Section\,\ref{sec:ZS_operators}.
For any $m\in\Z$, ${\mathcal D}_m$ is the closure of the domain bounded by $\Ga_m$.\footnote{Note that for
$|m|\ge N_0+1$, ${\mathcal D}_m$ is the disk $D_m$.}
\begin{Lemma}\hspace{-2mm}{\bf .}\label{Lem:auxiliary3}
Let $N\ge N_0$ and $R>0$. Then there exist an open neighborhood $\V$ of $\varphi_*$ in $\W$ and $0<C<\infty$ such that
for any $m\in\Z$ and for any $(\la,\si,\varphi)\in{\mathcal D}_m\times B_R\times\V$,
\[
\left|\prod_{|j|\ge N+1,j\ne m}\frac{\tsi_j-\la}{\sqrt[s]{(\la-\la_j^-(\varphi))(\la-\la_j^+(\varphi))}}\right|\le C\,,
\]
where $\tsi_j=j\pi+\si_j$ and $B_R=\{\si\in l^2\,|\,\|\si\|\le R\}$.
\end{Lemma}
\noindent{\em Proof.}
Consider first the case $|m|\ge N+1$.
Then for any $\si\in l^2$, $\varphi\in\W$, $\la\in D_m$, and $|j|\ge N+1$ with $j\ne m$
\begin{equation}\label{e:formula}
\frac{\tsi_j-\la}{\sqrt[s]{(\la-\la_j^-)(\la-\la_j^+)}}=
\frac{\tsi_j-\la}{\tau_j-\la}\cdot\frac{\tau_j-\la}{\sqrt[s]{(\la-\la_j^-)(\la-\la_j^+)}}\,,
\end{equation}
\begin{equation}\label{e:1+err1}
\frac{\tsi_j-\la}{\tau_j-\la}=1+\frac{\tsi_j-\tau_j}{\tau_j-\la}=1+O\Big(\frac{\tsi_j-\tau_j}{j-m}\Big)
\end{equation}
and by the choice of $\W$ in Section \ref{sec:ZS_operators}
\begin{eqnarray}\label{e:1+err2}
\frac{\tau_j-\la}{\sqrt[s]{(\la-\la_j^-)(\la-\la_j^+)}}=
\sqrt[+]{1+\frac{\ga_j^2}{4(\la-\la_j^-)(\la-\la_j^+)}}=
1+O\Big(\frac{\ga_j^2}{(m-j)^2}\Big)
\end{eqnarray}
where the constants in \eqref{e:1+err1} and \eqref{e:1+err2} depend only on the choice of the open neighborhood
$\W$ of $\varphi_*$ and $N$. 
Combining \eqref{e:formula} with \eqref{e:1+err1} and \eqref{e:1+err2} we get
\begin{eqnarray}
\frac{\tsi_j-\la}{\sqrt[s]{(\la-\la_j^-)(\la-\la_j^+)}}
&=&1+O\Big(\frac{|\si_j|+|j\pi-\tau_j|}{|m-j|}\Big)+O\Big(\frac{|\ga_j|^2}{|m-j|^2}\Big)\label{e:formula*1}\\
&=&1+O\Big(\frac{|\si_j|+|j\pi-\tau_j|+|\ga_j|}{|m-j|}\Big)\,.\label{e:formula*2}
\end{eqnarray}
Using Proposition \ref{Prop:periodic_spectrum} we choose a neighborhood $\V$ of
$\varphi_*$ in $\W$ so that for any $\varphi\in\V$, the sequence $(\la_j^{\pm}(\varphi)-j\pi)_{|j|\ge N_0+1}$ is bounded in $l^2$.
The statement of the Lemma for $|m|\ge N+1$ then follows from \eqref{e:formula*2} and inequality \eqref{e:r_m-rough}
of Lemma \ref{Lem:auxiliary2}. The case $|m|\le N$ is treated similarly.
\finishproof

\begin{Lemma}\hspace{-2mm}{\bf .}\label{Lem:auxiliary3'}
Let $N\ge N_0$ and $\si_0\in l^2$. Then for any $\ep>0$ there exist $m_0\ge 1$, an open neighborhood $U$ of $\si_0$
in $l^2$, and an open neighborhood $\V$ of $\varphi_*$ in $\W$ such that
for any $|m|\ge m_0$ and $(\la,\si,\varphi)\in{\mathcal D}_m\times U\times\V$
\[
\left|\,\,\,1-\prod_{|j|\ge N+1,j\ne m}\frac{\tsi_j-\la}{\sqrt[s]{(\la-\la_j^-(\varphi))(\la-\la_j^+(\varphi))}}\right|\le\ep\,,
\]
where $\tsi_j=j\pi+\si_j$.
\end{Lemma}
\noindent{\em Proof.} 
Arguing as in the proof of Lemma \ref{Lem:auxiliary3} we see that formula \eqref{e:formula*1} holds
for any $|m|\ge N_0+1$, $(\la,\si,\varphi)\in D_m\times l^2\times \W$, and $|j|\ge N+1$ with $j\ne m$. As a consequence,
\begin{equation}\label{e:formula*'}
\frac{\tsi_j-\la}{\sqrt[s]{(\la-\la_j^-)(\la-\la_j^+)}}=
1+O\Big(\frac{|\si_j|+|\tau_j-j\pi|+|\ga_j|^2}{|m-j|}\Big)
\end{equation}
where the constants depend only on the choice of $\W$ and $N$.
Finally, the Lemma follows from \eqref{e:formula*'}, inequality \eqref{e:r_m} of Lemma \ref{Lem:auxiliary2}, and
the continuity of the maps considered in Lemma \ref{Lem:some_analytic_maps} below.
\finishproof

The following Lemma is used in the proof of Lemma \ref{Lem:auxiliary3'}.
\begin{Lemma}\hspace{-2mm}{\bf .}\label{Lem:some_analytic_maps}
The maps defined on $\W$, 
\[
\varphi\mapsto(\tau_j(\varphi)-j\pi)_{|j|\ge N_0+1}\,\,\,\,\,\,\,\mbox{and}
\,\,\,\,\,\,\,\,\varphi\mapsto(\ga_j^2(\varphi))_{|j|\ge N_0+1}
\]
take values in $l^2$ and as such are analytic.
\end{Lemma}
\noindent{\em Proof.}
First, note that for any $|j|\ge N_0+1$, the mapping $\W\to\C$, $\varphi\mapsto\ga_j^2(\varphi)$,
is analytic. As
\[
\sum_{|j|\ge N_0+1}|\ga_j^2|^2\le\Big(\sum_{|j|\ge N_0+1}|\ga_j^2|\Big)^2
\]
we obtain from Proposition \ref{Prop:periodic_spectrum} that the map $\varphi\mapsto(\ga_j^2(\varphi))_{|j|\ge N_0+1}$
is locally bounded and hence analytic (see Theorem A.5 in \cite{KP}).
The analyticity of the map $\varphi\mapsto(\tau_j(\varphi)-j\pi)_{|j|\ge N_0+1}$ is proved in a similar way.
\finishproof

For $a,b\in\C$ let $[a,b]=\{ta+(1-t)b\,|\,t\in[0,1]\}$ and let $\Ga$ be a $C^1$-smooth, simple, closed curve
in $\C\setminus[a,b]$. The proof of the following lemma is straightforward and we omit it.
\begin{Lemma}\hspace{-2mm}{\bf .}\label{Lem:auxiliary4}
Let $f$ be holomorphic in an open neighborhood of $[a,b]$ containing $\Ga$.
Then
\[
\left|\oint_\Ga\frac{f(\la)}{\sqrt{(\la-a)(\la-b)}}\,d\la \right|\le 2\pi\max\limits_{\la\in[a,b]}|f(\la)|\,.
\]
\end{Lemma}

\section{Appendix B: Period map}\label{sec:periods}
In this Appendix we state and prove a result on the periods of a family of meromorphic differentials on
a compact Riemann surface which will be applied to prove Theorem \ref{Th:main}.
As this result is used at several instances we state it in a general form.
Let $\Sigma$ be a compact Riemann surface of genus $g\ge 0$ with the following data.

\vspace{0.3cm}

{\bf $(D1)$} $C_0,...,C_g$ are $C^1$-smooth, simple closed curves on $\Sigma$ 
dividing it into two connected components $\Sigma^\pm$,
\[
\Sigma^{-}\sqcup\Sigma^+=\Sigma\setminus(\sqcup_{j=0}^gC_k)
\]
which have the property that $C_k\cap C_l=\emptyset$ for any $0\le k,l\le g$, $k\ne l$.
Note that $C_0,...,C_g$ are the boundary cycles of $\Sigma^+$ in $\Sigma$.
The orientation of $C_0,...,C_g$ is chosen to be the one induced from the orientation on $\Sigma^+$ given by
the complex structure so that Stokes' formula holds on $\Sigma^+$.

\vspace{0.3cm}

{\bf $(D2)$} $\imath : \Sigma\to\Sigma$ is a holomorphic involution with $\imath : \Sigma^{\pm}\to\Sigma^{\mp}$ and
$\imath(C_k)=C_k$ for any $0\le k\le g$ as sets in $\Sigma$.

\vspace{0.3cm}

{\bf $(D3)$} $P_1^+,...,P_s^+$ are pairwise different points on $\Sigma^+$
and $P_j^-:=\imath(P^+_j)$ for any $1\le j\le s$. 

\vspace{0.3cm}

{\bf $(D4)$} For any $1\le k\le s$, $C_{g+k}$ is a simple closed $C^1$-smooth curve around $P_k^+$
which bounds a (small) open disk $U_{g+k}^+$ containing $P_k^+$ so that the closed disks
$\overline{U_{g+k}^+}$ ($1\le k\le s$) are contained in $\Sigma^+$ and pairwise disjoint.
The orientation of $C_{g+k}$ is induced from the orientation of $U_{g+k}^+$ so that Stokes' formula holds on
$U_{g+k}^+$.

\vspace{0.3cm}

{\bf $(D5)$} $\eta_1,...,\eta_{g+s}$ are {\em linearly independent} meromorphic differentials on $\Sigma$ so that
for any $1\le k\le g+s$
\begin{equation}\label{e:pull_back}
\imath^*(\eta_k)=-\eta_k\,\,\,\,\,\mbox{and}\,\,\,\,(\eta_k)\ge -D
\end{equation}
where $D$ is the divisor on $\Sigma$ given by
\[
\sum\limits_{j=1}^s (P^-_j+P^+_j).
\]

\vspace{0.3cm}

In particular, $(D5)$ means that the differentials $\eta_k$, $1\le k\le g+s$, may have poles only at
the points $P_{g+j}^\pm$, $1\le j\le s$, and that all these poles (if any) are of first order.

\vspace{0.3cm}

\begin{Lemma}\hspace{-2mm}{\bf .}\label{Lem:a-periods}
Assume that $(D1)-(D5)$ hold. Then for any $0\le k_0\le g+s$, the $(g+s)\times(g+s)$-matrix 
\[
X_{k_0}=\Big(\oint_{C_m}\eta_j\Big)_{1\le j\le g+s,0\le m\le g+s, m\ne k_0}
\]
is non-degenerate.
\end{Lemma}
{\em Proof.}
Assume that there exists $1\le k_0\le g+s$ so that the matrix $X_{k_0}$ is degenerate.
Then there exists a non-trivial linear combination,
\[
\eta:=\sum\limits_{k=0}^{g+s}c_k\eta_k\,\,\,\,\mbox{with}\,\,\,\,c_k\in\C
\]
so that
\begin{equation}\label{e:zero_periods}
\oint_{C_m}\eta=0\,\,\,\,\,\,\,\,\forall\,0\le m\le g+s,\,\,\,m\ne k_0\,.
\end{equation}
As $(\eta)\ge -D$, the Abelian differential $\eta$ is holomorphic in the interior of
$\Sigma^+\setminus(\sqcup_{k=1}^s U_{g+k}^+)$ and continuous on its boundary $\sqcup_{m=0}^{g+s} C_m$.
By Stokes' formula,
\[
\sum\limits_{m=0}^g\oint_{C_m}\eta-\sum\limits_{j=1}^s\oint_{C_{g+j}}\eta=0\,.
\]
Combined with \eqref{e:zero_periods} it then follows that 
\[
\oint_{C_{k_0}}\eta=0\,.
\] 
This, together with \eqref{e:zero_periods} and $(\eta)\ge-D$, implies that
$\eta$ is holomorphic on $\Sigma^+$. Taking into account that by $(D5)$, $\imath^*(\eta_k)=-\eta_k$ for any $1\le k\le g+s$,
it then follows that $\eta$ is a holomorphic differential on $\Sigma$.
As $\oint_{C_m}\eta=0$ for $1\le m\le g$ one then concludes that $\eta=0$, contradicting the linear independence of $\eta_k$'s
assumed in $(D5)$.
\finishproof

\end{document}